\documentclass[12pt,a4paper]{article}
\usepackage{amsmath}
\usepackage{amssymb}
\usepackage{cases}
\usepackage{bbm}
\usepackage{mathrsfs}
\usepackage{graphicx}
\usepackage{amsfonts}

\usepackage{verbatim}
\usepackage{enumerate}
\usepackage{theorem}
\usepackage{color}
\usepackage{cite}
\usepackage[pdfstartview=FitH]{hyperref}
\makeatletter
\def\tank#1{\protected@xdef\@thanks{\@thanks
 \protect\footnotetext[0]{#1}}}
\def\bigfoot{

 \@footnotetext}
\makeatother

\newcommand{\ea}{\end{array}}

\allowdisplaybreaks

\newtheorem{theorem}{Theorem}[section]
\newtheorem{lem}{Lemma}[section]
\newtheorem{prp}[theorem]{Proposition}
\newtheorem{thm}[theorem]{Theorem}
\newtheorem{cor}[theorem]{Corollary}
\newtheorem{dfn}[theorem]{Definition}

\newtheorem{remark}{Remark}[section]

\def\beq{\begin{equation}}
\def\nneq{\end{equation}}

\def\bthm{\begin{thm}}
\def\nthm{\end{thm}}

\def\blem{\begin{lem}}
\def\nlem{\end{lem}}
\def\bprf{\begin{proof}}
\def\nprf{\end{proof}}
\def\bprop{\begin{prop}}
\def\nprop{\end{prop}}
\def\brmk{\begin{rem}}
\def\nrmk{\end{rem}}

\def\bexa{\begin{exa}}
\def\nexa{\end{exa}}
\def\bcor{\begin{cor}}
\def\ncor{\end{cor}}

{\theorembodyfont{\rmfamily}
}

\title{Mckean-Vlasov stochastic differential equations with oblique reflection on non-smooth time dependent domains}

\footnotesize{\author{
Rong Wei$^{1,}$\thanks{wrbeauty@mail.ustc.edu.cn},\ \
 Saisai Yang$^{1,}$\thanks{yangss@mail.ustc.edu.cn},\ \
 Jianliang Zhai$^{1,}$\thanks{Zhaijl@ustc.edu.cn}\\
 {\em $^1$ School of Mathematical Sciences,}\\
 {\em University of Science and Technology of China,}\\
 {\em Hefei, 230026, China}
}
\date{}
\newenvironment{proof}{\par\noindent{\bf Proof:}}{\hspace*{\fill}$\blacksquare$\par}
\begin{document}
\maketitle
\noindent \textbf{Abstract:}
In this paper, we consider a class of Mckean-Vlasov stochastic differential equation with oblique reflection over an non-smooth time dependent domain. We establish the existence and uniqueness results of this class, address the propagation of chaos and prove a Fredlin-Wentzell type large deviations principle (LDP). One of the main difficulties is raised by the setting of non-smooth time dependent domain. To prove the LDP, a sufficient condition for the weak convergence method, which is suitable for Mckean-Vlasov stochastic differential equation,  plays an important role.


\vspace{4mm}


\vspace{3mm}
\numberwithin{equation}{section}
\section{Introduction}

Let $\mathcal{D}^{'}$ be a bounded connected domain in $\mathbb{R}^{1+d}$. For a given $T>0$, we define $$\mathcal{D}:=\mathcal{D}^{'}\cap([0,T]\times \mathbb{R}^{d})$$
as a time dependent domain. For any $t\in[0,T]$, we denote $\mathcal{D}_{t}:=\{x:(t,x)\in\mathcal{D}\}$ as the time sections of $\mathcal{D}$ and denote $\partial \mathcal{D}_{t}$ as the boundary of $\mathcal{D}_{t}$. Let $\overline{\mathcal{D}}_{t}$ be the closure of $\mathcal{D}_{t}$.
In this paper, we study a class of Mckean-Vlasov stochastic differential equation (SDE) with oblique reflection over the non-smooth time dependent domain $\mathcal{D}$:
\begin{align} \label{I-eq1.2}
\left\{
\begin{aligned}
& X_{t}=X_0+\int^{t}_{0}b(s,X_{s},\mu^X_{s})\mathrm{d}s+\int^{t}_{0}\sigma(s,X_{s},\mu^X_{s})\mathrm{d}W_{s}+K_{t},\ t\in[0,T],\\
& |K|_{t}=\int^{t}_{0}\mathbf{1}_{\{X_{s}\in\partial \mathcal{D}_{s}\}}\mathrm{d}|K|_{s},\quad K_{t}=\int^{t}_{0}\gamma(s,X_{s})\mathrm{d}|K|_{s}, \\
\end{aligned}
\right.
\end{align}
where the initial data $ X_0\in\overline{\mathcal{D}}_{0}$, $\mu^X_{s}$ is the law of $X_{s}$, $W=\{W_{t}\}_{t\in[0,T]}$ is an $m$-dimensional standard Wiener process defined on a complete probability space $(\Omega,\mathcal{F},\mathbb{F},\mathbb{P})$ with the filtration $\mathbb{F}:=\{\mathcal{F}_{t}\}_{t\in[0,T]}$ satisfying the usual conditions,  $\gamma$ is the so-called direction of reflection, and $|K|_t$ stands for the total variation of $K$ on $[0,t]$. For the precise conditions on $\mathcal{D}$, $b$, $\sigma$, and $\gamma$, we refer the reader to Section 2.

\vskip 0.1cm

Reflected SDEs have long been of interest of stochastic analysis since the works by Skorokhod \cite{Skorohod1961,Skorohod1962}.
There is an extensive literature on reflected SDEs over time independent domains; see, e.g., \cite{Lion,SaishoY,Costantini1990,DupuisIshii,Saisho} and the references therein.
The present paper is concerned with the setting of non-smooth time dependent domain, also known as ``non-cylindrical domain". The study of this topic is partly motivated by that reflecting Brownian motions in time dependent domains
arise in queueing theory \cite{[36],[43]}, statistical physics \cite{[13],[59]}, control theory \cite{[27],[28]}, and finance \cite{[26]}. There are only a few papers on this topic. We refer to \cite{Costantini2006,Nystr"om,LO,BCS2004,BurdzyW} for the related works.

Mckean-Vlasov/distribution-dependent SDEs with reflection over time
independent domains have also been studied. Motivated by non-linear Fokker-Planck equation with a Neumann boundary condition, the first paper to
study reflected McKean-Vlasov SDEs was \cite{SznitamnA}, proving existence, uniqueness and propagation of chaos.
The authors  did these in the general settings of a smooth bounded domain and bounded Lipschitz coefficients. In \cite{AdamsReisRavaille}, the authors established existence and uniqueness results for McKean-Vlasov SDEs constrained
to a convex domain  with coefficients that have superlinear growth in space and are non-Lipschitz in
measure.  They also proved propagation of chaos, a Freidlin-Wentzell type large deviation principle (LDP), and an Eyring-Kramer's law for the exit time from subdomains contained in the interior of the reflecting domain. In \cite{CDFM}, the authors applied a pathwise appproach to study existence, uniqueness and propagation of chaos for a class of reflected McKean-Vlasov SDEs with additive noise.
In \cite{GM1989}, the authors obtained existence, uniqueness and propagation of chaos for a class of Mckean-Vlasov SDEs with so-called sticky reflection. We also refer to \cite{HL,BEH, BCG} for reflected McKean-Vlasov SDEs with different kind of reflections.

To our knowledge, there are few results on reflected McKean-Vlasov SDEs in the setting of non-smooth time dependent domain.
In this setting, the aim of the present paper is to study a class of reflected Mckean-Vlasov SDE with oblique reflection. We first establish the existence and uniqueness of strong solution to (\ref{I-eq1.2}). We then  address the propagation of chaos, that is, we prove that the limit of a single
equation within the system of interacting equations (see (\ref{II-eq17}) below) converges to the dynamics of equation (\ref{I-eq1.2}). Lastly, we
prove a Freidlin-Wentzell's LDP for the strong solution to (\ref{I-eq1.2}). Our result on LDP is new even for the classical (i.e., distribution independent) reflected SDEs over time dependent domains.

Compared with the proof of the previous corresponding results for time independent domains, new difficulties occur, naturally, due to the fact that we are considering the setting of oblique reflection and non-smooth time dependent domain.  Sophisticated tools are needed. To obtain our results, we must employ involved  test functions $\{f_\delta\}_{\delta>0}$ (see Lemma \ref{lemII-2.1}) to deal with the reflection object rather than  $f_\delta(x)\equiv\|x\|^{p}$ used in the previous corresponding literature  for time independent domains. In particular, the approach used for the study of  Freidlin-Wentzell's LDP here is completely different from that in \cite{AdamsReisRavaille}.

Freidlin-Wentzell's LDP plays an important role in stochastic analysis, which can provide an exponential estimate of convergence for the probability of rare events
as the noise terms in stochastic systems tend to zero. Loosely speaking, it seeks a deterministic path such that the diffusion can be seen as a small random perturbation of this path. This type of LDP for classical/distribution-independent stochastic evolution equations has been extensively investigated;  see e.g., \cite{BD2019,DZ} and references therein. Large deviation problems for classical/distribution-independent reflected SDEs have been studied by several authors; see \cite{MST,WZZ2022}.

There are only a few results concerning the LDP for Mckean-Vlasov SDEs. The papers \cite{[RST],[HIP],[SY]} considered the case of the Gaussian
driving noise. Their proofs are based on exponential equivalence arguments. Some extra regularity
with respect to time on the coefficients is required. To obtain the LDP for  reflected Mckean-Vlasov SDE with normal reflection on time independent convex domain, the authors in \cite{AdamsReisRavaille} adopted the exponential equivalence arguments, certain time discretization and approximating technique, assuming that the coefficients satisfy some
extra time $\rm H\ddot{o}lder$ continuity conditions; see Assumption 4.1 in \cite{AdamsReisRavaille}. These approaches are very difficult to be applied to the case of infinite dimensional situations and the case of L\'evy driving noise, and requires stronger conditions on the coefficients as mentioned above. The use of these approaches  to deal with our problem would be very hard, if not impossible at all.

In a recent paper \cite{LSZ}, the third author and his collaborators presented a sufficient condition  to
verify the criteria of the weak convergence method \cite{BD2019}. The sufficient condition is suitable for Mckean-Vlasov SDE in finite and infinite dimensions. They then applied the sufficient condition to establish an LDP for distribution-dependent SDE driven by L\'evy process. It is the first paper to fully use the weak convergence method to establish distribution-dependent SDE, and it requires the very natural Lipschitz conditions on the coefficients without the extra assumptions appearing in the literatures mentioned above.  We would like to point out the article \cite{[HLL2021]} in which the authors already applied the sufficient condition introduced in \cite{LSZ}  to infinite dimensional situations. In this paper, we apply the sufficient condition introduced in \cite{LSZ} to study Freidlin-Wentzell's LDP in our setting. Our proof seems to be smoother than that in \cite{AdamsReisRavaille}, and no extra regularity with respect to time on the coefficients is required.

Finally, we point out that the weak convergence method  is proved to be a powerful tool to establish LDPs for various dynamical systems driven by Gaussian noise and/or Poisson random measures. A listing of  the applications can be found
 in \cite{LSZ}.

\vskip 0.3cm
An outline of the current work is as follows. In Section 2, we  prepare some basic concepts, notations and assumptions. The existence and uniqueness of the equation (\ref{I-eq1.2}) are established in Section 3. After that, we give the result of the propagation of chaos in Section 4. Finally, in Section 5, we study the LDP for (\ref{I-eq1.2}).


\section{Preliminary and Assumptions}\label{secII}
Throughout this paper we will use the following notations and assumptions.

\vskip 0.3cm

For $x,y\in\mathbb{R}^{d}$, let $\langle x,y\rangle$ and $\|x\|:=\langle x,x\rangle^{\frac{1}{2}}$ stand for the Euclidean inner product and norm, respectively. For any $a\in \mathbb{R}^{d}$ and $r>0$, we denote $B(a,r):=\{x\in\mathbb{R}^{d}:\|x-a\|\leq r\}$, and $S(a,r):=\{x\in\mathbb{R}^{d}:\|x-a\|=r\}$. For any metric spaces $E$ and $F$, we define the following spaces of functions mapping $E$ (or $[0,T]\times E$) to $F$. $C(E,F)$ denotes the set of continuous functions. Let $C^{1,2}([0,T]\times E,F)$ denote the
set of functions, whose elements are
continuously differentiable once with respect to the time variable and twice with respect to any
space variable, and by $C^{1,2}_b([0,T]\times E,F)$ we denote the space of bounded functions in $C^{1,2}([0,T]\times E,F)$
having bounded derivatives.
Denote by $\mathcal{W}^{1,p}([0,T],[0,\infty))$ the Sobolev space of functions whose first order weak derivatives belong to $L^p([0,T],\mathbb{R})$, here $L^p([0,T],\mathbb{R})$ is the $\mathbb{R}$-valued $L^p$-space.
Analogously, we define $C^{1,2,2}\left([0,T]\times\mathbb{R}^{d}\times\mathbb{R}^{d},\mathbb{R}\right)$, $C^{1,2}(\overline{\mathcal{D}},\mathbb{R})$, $C^{1,2}_{b}(\mathbb{R}^{1+d},B(0,1))$, etc. Here $\overline{\mathcal{D}}$
is closure of $\mathcal{D}$.
For $u\in C^{1,2,2}\left([0,T]\times\mathbb{R}^{d}\times\mathbb{R}^{d},\mathbb{R}\right)$, $D_t u$ denotes the differential with respect to the time variable, $D_x u$ and $D_y u$ denote the gradient of $u$ with respect to the spatial variables $x$ and $y$, respectively, and $D^2u$ denote the Hessian matrix of $u$ with respect to spatial variables.
And we use the letter $C$ for a generic positive constant whose value may change from line to line.

\vskip 0.3cm

Let $(E,m)$ be a metric space with Borel $\sigma$-algebra $\mathcal{B}(E)$. We use $\mathcal{P}^{2}(E)$ to denote the set of all probability measures on $(E,\mathcal{B}(E))$ which have finite moment of order 2. It is well known, see, e.g., \cite{Bolley,Villani}, that $\mathcal{P}^{2}(E)$ is a complete metric space under the Wasserstein 2-distance
\begin{equation}
\mathcal{W}_{E}(\mu,\nu)=\left\{\inf_{\pi\in\Pi(\mu,\nu)}\left[\int_{E\times E}
m(x,y)^{2}\pi(\mathrm{d}x,\mathrm{d}y)\right]\right\}^{\frac{1}{2}},\ \mu,\nu\in \mathcal{P}^{2}(E).
\end{equation}
Here and throughout this paper $\Pi(\mu,\nu)\subset \mathcal{P}^{2}(E\times E)$ is the set of all joint distributions over $E\times E$ with marginals  $\mu$ and $\nu$.

Denote by $\mathcal{X}$ the space of $\mathbb{R}^{d}$-valued continuous functions $f$ on $[0,T]$ satisfying $f(l)\in\overline{\mathcal{D}}_{l}$ for each $l\in[0,T]$, and $\mathcal{X}$ is equipped with the uniform topology, then it is a Polish space.
The symbols $\mathcal{P}^{2}(\overline{\mathcal{D}}_{t})$, $\mathcal{W}_{\overline{\mathcal{D}}_{t}}(\cdot,\cdot)$, $\mathcal{P}^{2}(\mathbb{R}^d)$, $\mathcal{W}_{\mathbb{R}^d}(\cdot,\cdot)$,  $\mathcal{P}^{2}(\mathcal{X})$, and $\mathcal{W}_{\mathcal{X}}(\cdot,\cdot)$ will be used in the sequel.

For any $\mu\in\mathcal{P}^{2}(\mathcal{X})$,  there is a natural surjection
                $$\mathcal{P}^{2}(\mathcal{X})
                \ni\mu\mapsto(\mu_{t})_{t\in[0,T]}\in C([0,T];
                \mathcal{P}^{2}(\mathbb{R}^d)),$$
                where $\mu_{t}$ is the pushforward  measure with respect to path evaluation defined by
                $$\mu_{t}(A):=\int_{\mathcal{X}}\mathbf{1}_{\{x \in \mathcal{X};x(t)\in A\}}(x)\mu(\mathrm{d}x), \ \forall A \in \mathcal{B}(\mathbb{R}^d).$$
                Note that $\mu_{t}\in \mathcal{P}^{2}(\overline{\mathcal{D}}_{t})$, and it is easy to see that for any $\mu,\nu\in\mathcal{P}^{2}(\mathcal{X})$,
                \begin{equation}\label{III-eq3.1 00}
                \mathcal{W}_{\mathbb{R}^d}(\mu_{t},\nu_{t})= \mathcal{W}_{\overline{\mathcal{D}}_{t}}(\mu_{t},\nu_{t}),
                \end{equation}
                and

                \begin{equation}\label{III-eq3.1}
                \sup_{t\in[0,T]}\mathcal{W}_{\overline{\mathcal{D}}_{t}}(\mu_{t},\nu_{t})
                =
                \sup_{t\in[0,T]}\mathcal{W}_{\mathbb{R}^d}(\mu_{t},\nu_{t})
                \leq
                 \mathcal{W}_{\mathcal{X}}(\mu,\nu).
                \end{equation}

\vskip 0.3cm

At the end of this section, we give some assumptions about the boundary conditions and coefficients. Throughout this work we will assume that for any $t\in [0,T]$, $\mathcal{D}_{t}$ satisfies
\begin{equation*}
\mathcal{D}_{t}\neq \emptyset\ \mbox{and that}\ \mathcal{D}_{t}\ \mbox{is a bounded connected set for every}\ t\in [0,T],
\end{equation*}
and the direction of reflection at $x\in\partial{\mathcal{D}_{t}}$ given by $\gamma(t,x)$ satisfies
\begin{equation*}
\gamma\in C^{1,2}_{b}(\mathbb{R}^{1+d},B(0,1)),\label{II-eq2.2}
\end{equation*}
and
\begin{equation*}
\gamma(t,x)\in S(0,1)\ \mbox{for all}\ (t,x)\in V,
\end{equation*}
where $V$ is an open set containing $[0,T]\times\mathbb{R}^{d}\backslash \mathcal{D}$. In particular, there exists a constant $\rho\in (0,1)$ such that the exterior cone condition
\begin{equation}\label{II-eq2}
\bigcup_{0\leq\xi\leq \rho}B(x-\xi\gamma(t,x),\xi \rho)\subset\mathcal{D}^{c}_{t}, \ \mbox{for all}\ x\in\partial \mathcal{D}_{t},t\in[0,T]
\end{equation}
holds. Let $d(t,x):=\inf_{y\in\mathcal{D}_{t}}\|x-y\|$ for any $t\in [0,T]$, $x\in\mathbb{R}^{d}$ and assume that for some fixed $p\in(1,\infty)$ and all $x\in\mathbb{R}^{d}$, $d(\cdot,x)\in\mathcal{W}^{1,p}([0,T],[0,\infty))$ with Sobolev norm uniformly bounded in space. We also assume that the time derivative of $d(t,x)$ is jointly measurable in $(t,x)$.
\vskip 0.3cm
Now we give the assumption on the drift and diffusion coefficients.
Let $b:[0,T]\times\mathbb{R}^{d}\times\mathcal{P}^{2}(\mathbb{R}^{d})\to\mathbb{R}^{d}$ and $\sigma:[0,T]\times\mathbb{R}^{d}\times\mathcal{P}^{2}(\mathbb{R}^{d})\to \mathcal{L}_2(\mathbb{R}^m,\mathbb{R}^d)$ be measurable functions,
here $\mathcal{L}_2(\mathbb{R}^m,\mathbb{R}^d)=\mathbb{R}^{m}\otimes \mathbb{R}^d$ is the space of all Hilbert-Schmidt operators from $\mathbb{R}^m$ to $\mathbb{R}^d$ equipped with
the usual Hilbert-Schmidt norm $\|\cdot\|_{\mathcal{L}_2}$. When there is no danger of causing ambiguity, we denote $\|\cdot\|=\|\cdot\|_{\mathcal{L}_2}$.

Let us give the following assumptions about the coefficients in \eqref{I-eq1.2}.
\begin{itemize}
\item [(A1)] There exists a function $L\in L^{2}([0,T],\mathbb{R}^+)$ such that for any $t\in [0,T]$, $x,y\in \mathbb{R}^{d}$, and $\mu,\nu\in\mathcal{P}^{2}(\mathbb{R}^{d})$,
\begin{equation}
\|b(t,x,\mu)-b(t,y,\nu)\|
+
\|\sigma(t,x,\mu)-\sigma(t,y,\nu)\|
\leq
L(t)(\|x-y\|+\mathcal{W}_{\mathbb{R}^{d}}(\mu,\nu))\label{II-eq2.16}.
\end{equation}
\end{itemize}
\begin{remark}\label{Remark2.1}
The Assumption (A1) implies that for any $t\in[0,T]$, $x\in\overline{\mathcal{D}_{t}}$, $\mu\in\mathcal{P}^{2}(\overline{\mathcal{D}_{t}})$, the coefficients $b$ and $\sigma$ satisfy
\begin{equation}\label{Remeq2.1}
\|b(t,x,\mu)\|+\|\sigma(t,x,\mu)\|\leq C(1+L(t)).
\end{equation}
\end{remark}

Due to Theorem 2.7 in \cite{LO}, we have the following two results.

\begin{prp}\label{prop 1}
  Under Assumption (A1), for each fixed $\mu\in\mathcal{P}^{2}(\mathcal{X})$ and initial condition $X_0\in\overline{\mathcal{D}}_{0}$,
there exists a unique strong solution $X^{\mu}=\{X^{\mu}_t,t\in[0,T]\}$ to the following reflected SDE:
\begin{eqnarray}\label{eq Pop 1}
& & X^{\mu}_{t}=X_0+\int^{t}_{0}b(s,X^{\mu}_{s},\mu_{s})\mathrm{d}s+\int^{t}_{0}\sigma(s,X^{\mu}_{s},\mu_{s})\mathrm{d}W_{s}+K^{\mu}_{t},\ t\in[0,T],\nonumber
   \\
& & X^{\mu}_{t}\in\overline{\mathcal{D}}_{t},\quad |K^{\mu}|_{t}=\int^{t}_{0}{\bf 1}_{\{X^{\mu}_{s}\in\partial\mathcal{D}_{s}\}}\mathrm{d}|K^{\mu}|_{s}<\infty,\quad
    \ K^{\mu}_{t}=\int^{t}_{0}\gamma(s,X^{\mu}_{s})\mathrm{d}|K^{\mu}|_{s},
\end{eqnarray}
and $\mathbb{E}[\sup_{0\leq t\leq T}\|X^{\mu}_{t}\|^{2}]<\infty.$
\end{prp}

\begin{prp}\label{prop 2}
  Assume that  there exists a function $\widetilde{L}\in L^{1}([0,T],\mathbb{R}^+)$ such that for any $t\in [0,T]$, $x,y\in \mathbb{R}^{d}$, and $\nu\in\mathcal{P}^{2}(\mathbb{R}^{d})$,
\begin{equation}
\|b(t,x,\nu)-b(t,y,\nu)\|
\leq
\widetilde{L}(t)\|x-y\|.
\end{equation}
Then for each fixed $\mu\in\mathcal{P}^{2}(\mathcal{X})$ and initial condition $X_0\in\overline{\mathcal{D}}_{0}$,
there exists a unique strong solution $\widetilde{X}^{\mu}=\{\widetilde{X}^{\mu}_t,t\in[0,T]\}$ to the following reflected differential equation:
\begin{eqnarray*}
& & \widetilde{X}^{\mu}_{t}=X_0+\int^{t}_{0}b(s,\widetilde{X}^{\mu}_{s},\mu_{s})\mathrm{d}s+\widetilde{K}^{\mu}_{t},\ t\in[0,T],\nonumber
   \\
& & \widetilde{X}^{\mu}_{t}\in\overline{\mathcal{D}}_{t},\quad |\widetilde{K}^{\mu}|_{t}=\int^{t}_{0}\mathbf{1}_{\{\widetilde{X}^{\mu}_{s}\in\partial\mathcal{D}_{s}\}}\mathrm{d}|\widetilde{K}^{\mu}|_{s}<\infty,\quad
    \ \widetilde{K}^{\mu}_{t}=\int^{t}_{0}\gamma(s,\widetilde{X}^{\mu}_{s})\mathrm{d}|\widetilde{K}^{\mu}|_{s},
\end{eqnarray*}
and $\sup_{0\leq t\leq T}\|\widetilde{X}^{\mu}_{t}\|^{2}<\infty.$
\end{prp}

\begin{remark}
  The assumptions of the above two propositions are a little different with that of Theorem 2.7 in \cite{LO}, in which it requires that $L(t)\equiv\widetilde{L}(t)\equiv K$; see (2.17) in \cite{LO}. However, using similar arguments in the proof of Theorem 2.7 in \cite{LO} and the idea in
  the proof of Theorem \ref{II-thm2.2} in this paper, it is not difficult to get the above two propositions, and the proofs are omitted here.
\end{remark}

\section{Existence and Uniqueness for RMVSDE}

In this section, we establish the existence and uniqueness of strong solution to equation (\ref{I-eq1.2}).

We now introduce the definition of the solution to equation (\ref{I-eq1.2}).
\begin{dfn}
 An $\mathbb{R}^d$-valued stochastic  process $X=\{X_t,t\in[0,T]\}$ is called a strong solution to (\ref{I-eq1.2}) in $\overline{\mathcal{D}}$ driven by the Wiener process $W$ and with coefficients $b$ and $\sigma$, direction of reflection along $\gamma$ and initial condition $X_0\in\overline{\mathcal{D}}_{0}$, if $X$ is an $\mathbb{F}$-adapted stochastic process which satisfies, $\mathbb{P}$-a.s.,
 $$
 X\in\mathcal{X},
 $$
 and, for any $t\in[0,T]$,
 $$X_{t}=X_0+\int^{t}_{0}b(s,X_{s},\mu^X_{s})\mathrm{d}s+\int^{t}_{0}\sigma(s,X_{s},\mu^X_{s})\mathrm{d}W_{s}+K_{t},$$
 where $\mu^X_{s}$ is the law of $X_{s}$,
 \begin{equation}
  |K|_{t}=\int^{t}_{0}\mathbf{1}_{\{X_{s}\in\partial \mathcal{D}_{s}\}}\mathrm{d}|K|_{s}<\infty,\quad \text{and }\quad
 K_{t}=\int^{t}_{0}\gamma(s,X_{s})\mathrm{d}|K|_{s}.
 \end{equation}
\end{dfn}

Before giving the main result in this section, we first present the following lemma, which is inspired by Lemmas 3.2 and 3.3 in \cite{{LO}}. It plays an important role in this work.
\begin{lem}\label{lemII-2.1}
\begin{itemize}
\item[(1)]There exist positive constants $\chi$ and $C$ (independent of $\delta$) and a family of the so-called test functions $\{f_\delta\}_{\delta>0}\subseteq C^{1,2,2}\left([0,T]\times\mathbb{R}^{d}\times\mathbb{R}^{d},\mathbb{R}\right)$ such that for all $(t,x,y)\in[0,T]\times\mathbb{R}^{d}\times\mathbb{R}^{d}$,
    \begin{equation}\label{II-eq2.4}
    f_{\delta}(t,x,y)\geq \chi\frac{\|x-y\|^{2}}{\delta},
    \end{equation}
    \begin{equation}\label{II-eq2.5}
    f_{\delta}(t,x,y)\leq C(\delta+\frac{\|x-y\|^{2}}{\delta}),
    \end{equation}
    \begin{equation}\label{II-eq2.6}
    \langle D_{x}f_{\delta}(t,x,y),\gamma(t,x)\rangle\leq C\frac{\|x-y\|^{2}}{\delta},\ \ \mbox{for}\ x\in\partial\mathcal{D}_{t}, \ y\in\overline{D}_{t},
    \end{equation}
    \begin{equation}\label{II-eq2.8}
    \langle D_{y}f_{\delta}(t,x,y),\gamma(t,y)\rangle \leq C\frac{\|x-y\|^{2}}{\delta},\ \
    \mbox{for}\ y\in\partial\mathcal{D}_{t},\ \ x\in \overline{D}_{t},
    \end{equation}
    \begin{equation}\label{II-eq2.9}
    \|D_{t}f_{\delta}(t,x,y)\|\leq C\frac{\|x-y\|^{2}}{\delta},
    \end{equation}
    \begin{equation}\label{II-eq2.10}
    \|D_{y}f_{\delta}(t,x,y)\|\leq C\frac{\|x-y\|}{\delta},\ \
    \|D_{x}f_{\delta}(t,x,y)+D_{y}f_{\delta}(t,x,y)\|\leq C\frac{\|x-y\|^{2}}{\delta},
    \end{equation}
    \begin{equation}\label{II-eq2.11}
    D^{2}f_{\delta}(t,x,y)\leq\frac{C}{\delta}
    \left(\begin{array}{cc}
    I & -I \\
    -I & I
    \end{array}\right)
    +\frac{C\|x-y\|^{2}}{\delta}
    \left(\begin{array}{cc}
    I & 0 \\
    0 & I
    \end{array}\right).
    \end{equation}
        Here $I$ and $0$ denote the unit matrix and zero matrix of size $d\times d$, respectively,  and for any $d\times d$ real
symmetric matrices $X$ and $Y$, we write $X\leq Y$ if $\langle (X-Y)\xi,\xi\rangle\leq 0$ for all $\xi\in\mathbb{R}^d$.
\item[(2)] There exists a nonnegative function $\alpha\in C^{1,2}(\overline{\mathcal{D}},\mathbb{R})$, which satisfies
\begin{equation}\label{II-eq2.12}
\langle D_{x}\alpha(t,x),\gamma(t,x)\rangle\geq 1,\ \mbox{for any}\ x\in\partial\mathcal{D}_{t}, t\in[0,T].
\end{equation}
\end{itemize}
\end{lem}

\begin{proof}
By \eqref{II-eq2}, there exist positive constants $\theta$ and $\kappa$ such that for any $(t,y)\in[0,T]\times\overline{D}_{t}$, $x\in\partial D_{t}$ satisfying $\|x-y\|<\kappa$, we have
\begin{equation}\label{Kap}
\langle y-x,\gamma(t,x)\rangle\geq-\theta\|x-y\|.
\end{equation}
Let $\{f_{\delta}\}_{\delta>0}$ be constructed same as $\{w_{\epsilon}\}_{\epsilon>0}$ in Lemma 3.2 and $\alpha$ defined in Lemma 3.3 in \cite{{LO}}. Then by Lemmas 3.2 and 3.3 in \cite{{LO}}, we only need to give the proof for \eqref{II-eq2.6} and \eqref{II-eq2.8}.

For $(t,y)\in[0,T]\times\overline{D}_{t}$ and $x\in\partial D_{t}$, when $\|x-y\|<\kappa$, it is easy to see that \eqref{II-eq2.6} holds following from \eqref{Kap} and (3.16) in \cite{{LO}}. When $\|x-y\|\geq\kappa$, by \eqref{II-eq2.10} we have
\begin{equation*}
\langle D_{x}f_{\delta}(t,x,y),\gamma(t,x)\rangle \leq\frac{C\|x-y\|}{\delta}\leq \frac{C\|x-y\|^{2}}{\kappa\delta}.
\end{equation*}
Therefore, we obtain \eqref{II-eq2.6}.
Similarly, we can also get \eqref{II-eq2.8}.
\end{proof}

\begin{thm}\label{II-thm2.2}
Under Assumption (A1), there exists a unique strong solution to (\ref{I-eq1.2}).
\end{thm}

\begin{proof}
We will make use of the Banach fixed point theorem to prove this result.

Define a truncated Wasserstein distance on $\mathcal{P}^{2}(\mathcal{X})$ by
  \begin{equation}
                \mathcal{W}_{t}(\mu,\nu)=\left\{\inf_{\pi\in\Pi(\mu,\nu)}\left[\int_{\mathcal{X}\times \mathcal{X}}
                \sup_{s\in[0,t]}\|x_{s}-y_{s}\|^{2}\pi(\mathrm{d}x,\mathrm{d}y)\right]\right\}^{\frac{1}{2}}.
                \end{equation}
It is easy to see that $\mathcal{W}_{T}(\mu,\nu)=\mathcal{W}_{\mathcal{X}}(\mu,\nu)$.


By Proposition \ref{prop 1}, define a map $\Phi:\mathcal{P}^{2}(\mathcal{X})\to \mathcal{P}^{2}(\mathcal{X})$ by setting $\Phi(\mu)=Law(X^{\mu})$,
where $X^{\mu}$ is the unique solution to (\ref{eq Pop 1}) and $Law(X^{\mu})$ is the law of the solution $X^{\mu}$ on $\mathcal{X}$.  We give the existence and uniqueness of (\ref{I-eq1.2}) in $[0,T]$ by proving that the mapping $\Phi$ has a unique fixed point in $\mathcal{P}^{2}(\mathcal{X})$.
Define
\begin{equation}\label{V}
V(t,x,y):=\exp\{-\lambda (\alpha(t,x)+\alpha(t,y))\}f_{\delta}(t,x,y):=u(t,x,y)f_{\delta}(t,x,y),
\end{equation}
where the function $\alpha$ and $f_{\delta}$ are given in Lemma \ref{lemII-2.1}. For any $\mu,\nu\in\mathcal{P}^{2}(\mathcal{X})$, let
$X^\nu$ be the unique solution to (\ref{eq Pop 1}) with $\mu$ replaced by $\nu$. By It\^o's formula, we can make decomposition of $V(t,X^{\mu}_{t},X^{\nu}_{t})$ as follows
\begin{equation}\label{III-eq3.2}
V(t,X^{\mu}_{t},X^{\nu}_{t})=V(0,X^{\mu}_{0},X^{\nu}_{0})+\sum_{i=1}^{6}I_{i}(t),
\end{equation}
where
\begin{eqnarray}
I_{1}(t)&:=&    \int^{t}_{0}D_{s}V(s,X^{\mu}_{s},X^{\nu}_{s})ds,\nonumber\\
I_{2}(t)&:=&    \int^{t}_{0}\langle D_{x}u(s,X^{\mu}_{s},X^{\nu}_{s}),b(s,X^{\mu}_{s},\mu_{s})\rangle
                f_{\delta}(s,X^{\mu}_{s},X^{\nu}_{s})\mathrm{d}s\nonumber\\
        &  &    +\int^{t}_{0}\langle D_{y}u(s,X^{\mu}_{s},X^{\nu}_{s}),b(s,X^{\nu}_{s},\nu_{s})\rangle
                f_{\delta}(s,X^{\mu}_{s},X^{\nu}_{s})\mathrm{d}s,\nonumber\\
I_{3}(t)&:=&    \int^{t}_{0}\langle D_{x}f_{\delta}(s,X^{\mu}_{s},X^{\nu}_{s})+
                D_{y}f_{\delta}(s,X^{\mu}_{s},X^{\nu}_{s}),b(s,X^{\mu}_{s},\mu_{s})\rangle
                u(s,X^{\mu}_{s},X^{\nu}_{s})\mathrm{d}s\nonumber\\
        &  &    +\int^{t}_{0}\langle D_{y}f_{\delta}(s,X^{\mu}_{s},X^{\nu}_{s}),b(s,X^{\nu}_{s},\nu_{s})-b(s,X^{\mu}_{s},\mu_{s})\rangle
                u(s,X^{\mu}_{s},X^{\nu}_{s})\mathrm{d}s\nonumber\\
        I_{4}(t)&:=&    \frac{1}{2}\int^{t}_{0}tr \left[\left(
                \begin{array}{cc}
                \sigma(s,X^{\mu}_{s},\mu_{s}) \\
                \sigma(s,X^{\nu}_{s},\nu_{s}) \\
                \end{array}
                \right)^{T}D^{2}V(s,X^{\mu}_{s},X^{\nu}_{s})
                \left(
                \begin{array}{cc}
                \sigma(s,X^{\mu}_{s},\mu_{s}) \\
                \sigma(s,X^{\nu}_{s},\nu_{s}) \\
                \end{array}
                \right)\right]\mathrm{d}s,\nonumber
\end{eqnarray}
where $A^{T}$ denotes the transpose of a matrix $A$,
\begin{eqnarray}
I_{5}(t)&:=&    \int^{t}_{0}\langle D_{x}f_{\delta}(s,X^{\mu}_{s},X^{\nu}_{s}),\gamma(s,X^{\mu}_{s})\rangle u(s,X^{\mu}_{s},X^{\nu}_{s})
                \mathrm{d}|K^{\mu}|_{s}\nonumber\\
        &  &    +\int^{t}_{0}\langle D_{x}u(s,X^{\mu}_{s},X^{\nu}_{s}),\gamma(s,X^{\mu}_{s})\rangle
                f_{\delta}(s,X^{\mu}_{s},X^{\nu}_{s})\mathrm{d}|K^{\mu}|_{s}\nonumber\\
        &  &    +\int^{t}_{0}\langle D_{y}f_{\delta}(s,X^{\mu}_{s},X^{\nu}_{s}),\gamma(s,X^{\nu}_{s})\rangle
                u(s,X^{\mu}_{s},X^{\nu}_{s})\mathrm{d}|K^{\nu}|_{s}\nonumber\\
        &  &    +\int^{t}_{0}\langle D_{y}u(s,X^{\mu}_{s},X^{\nu}_{s}),\gamma(s,X^{\nu}_{s})\rangle
                f_{\delta}(s,X^{\mu}_{s},X^{\nu}_{s})\mathrm{d}|K^{\nu}|_{s},\nonumber\\
I_{6}(t)&:=&    \int^{t}_{0}\langle D_{x}u(s,X^{\mu}_{s},X^{\nu}_{s}),\sigma(s,X^{\mu}_{s},\mu_{s})\rangle
                f_{\delta}(s,X^{\mu}_{s},X^{\nu}_{s})\mathrm{d}W_{s}\nonumber\\
        &  &    +\int^{t}_{0}\langle D_{y}u(s,X^{\mu}_{s},X^{\nu}_{s}),\sigma(s,X^{\nu}_{s},\nu_{s})\rangle
                f_{\delta}(s,X^{\mu}_{s},X^{\nu}_{s})\mathrm{d}W_{s}\nonumber\\
        &  &    +\int^{t}_{0}\langle D_{x}f_{\delta}(s,X^{\mu}_{s},X^{\nu}_{s})+D_{y}f_{\delta}(s,X^{\mu}_{s},X^{\nu}_{s}),
                \sigma(s,X^{\mu}_{s},\mu_{s})\rangle u(s,X^{\mu}_{s},X^{\nu}_{s})\mathrm{d}W_{s}\nonumber\\
        &  &    +\int^{t}_{0}\langle D_{y}f_{\delta}(s,X^{\mu}_{s},X^{\nu}_{s}),\sigma(s,X^{\nu}_{s},\nu_{s})
                -\sigma(s,X^{\mu}_{s},\mu_{s})\rangle
                u(s,X^{\mu}_{s},X^{\nu}_{s})\mathrm{d}W_{s}.\nonumber
\end{eqnarray}

Now, we are going to estimate $I_{1}(t)-I_{6}(t)$ one by one.
From (\ref{II-eq2.5}), (\ref{II-eq2.9}) and the regularity of $u$, one can see that
\begin{equation}
   I_{1}(t)
\leq C(\lambda)\int^{t}_{0}(\delta+\frac{\|X^{\mu}_{s}-X^{\nu}_{s}\|^{2}}{\delta})\mathrm{d}s\label{III-eq3.4}.
\end{equation}
By the regularity of $u$ and (\ref{II-eq2.5}),
\begin{equation}\label{III-eq3.2'}
I_{2}(t)\leq C(\lambda)\int^{t}_{0}(1+L(s))(\delta+\frac{\|X^{\mu}_{s}-X^{\nu}_{s}\|^{2}}{\delta})\mathrm{d}s.
\end{equation}
Based on (\ref{II-eq2.10}) and (\ref{II-eq2.16}),
\begin{eqnarray}
I_{3}(t)
&\leq&  C(\lambda)\int^{t}_{0}(1+L(s))\frac{\|X^{\mu}_{s}-X^{\nu}_{s}\|^{2}}{\delta}\mathrm{d}s\nonumber\\
& &     +C(\lambda)\int^{t}_{0}L(s)
        \frac{\|X^{\mu}_{s}-X^{\nu}_{s}\|}{\delta}(\|X^{\mu}_{s}-X^{\nu}_{s}\|+\mathcal{W}_{\overline{\mathcal{D}}_{s}}(\mu_{s},\nu_{s}))\mathrm{d}s\\
&\leq&  C(\lambda)\int_{0}^{t}(1+L(s))(\frac{\|X^{\mu}_{s}-X^{\nu}_{s}\|^{2}}{\delta}
        +\frac{\mathcal{W}_{\overline{\mathcal{D}}_{s}}(\mu_{s},\nu_{s})^{2}}{\delta})\mathrm{d}s.\nonumber
\end{eqnarray}
From (\ref{II-eq2.11}), there is a constant $C(\lambda)>0$ such that
\begin{equation*}
D^{2}V(t,x,y)\leq C(\lambda)\left[\frac{1}{\delta}
\left(
    \begin{array}{cc}
    I & -I \\
    -I & I \\
    \end{array}
\right)+
(\delta+\frac{\|x-y\|^{2}}{\delta})
\left(
    \begin{array}{cc}
    I & 0 \\
    0 & I \\
    \end{array}
\right)
\right].
\end{equation*}
Therefore
\begin{equation}
I_{4}(t)\leq C(\lambda)\int^{t}_{0}(1+L^{2}(s))\left[\delta+\frac{\|X^{\mu}_{s}-X^{\nu}_{s}\|^{2}}{\delta}
        +\frac{\mathcal{W}_{\overline{\mathcal{D}}_{s}}(\mu_{s},\nu_{s})^{2}}{\delta}\right]\mathrm{d}s\label{III-eq3.8}.
\end{equation}
Due to (\ref{II-eq2.4}), (\ref{II-eq2.6}), (\ref{II-eq2.8}) and (\ref{II-eq2.12}),
\begin{equation}
I_{5}(t)\leq(C-\lambda\chi)\int^{t}_{0}\frac{\|X^{\mu}_{s}-X^{\nu}_{s}\|^{2}}{\delta}\mathrm{d}|K^{\mu}|_{s}+
(C-\lambda\chi)\int^{t}_{0}\frac{\|X^{\mu}_{s}-X^{\nu}_{s}\|^{2}}{\delta}\mathrm{d}|K^{\nu}|_{s}.
\end{equation}
So, by putting $\lambda=\frac{C}{\chi}$ all integral with respect to $|K^{\mu}|$ and $|K^{\nu}|$ vanish.
Dropping the $\lambda$ dependence from the constants, (\ref{II-eq2.4}) and (\ref{III-eq3.2})-(\ref{III-eq3.8}) imply
\begin{eqnarray}
 \frac{1}{C} \frac{\|X^{\mu}_{t}-X^{\nu}_{t}\|^{2}}{\delta}
&\leq&  V(t,X^{\mu}_{t},X^{\nu}_{t})\nonumber\\
&\leq& V(0,X^{\mu}_{0},X^{\nu}_{0})+C\delta\int^{t}_{0}(1+L^{2}(s))\mathrm{d}s\nonumber\\
& & +C\int^{t}_{0}(1+L^{2}(s))(\frac{\|X^{\mu}_{s}-X^{\nu}_{s}\|^{2}}{\delta}+\frac{\mathcal{W}_{\overline{\mathcal{D}}_{s}}(\mu_{s},\nu_{s})^{2}}{\delta})\mathrm{d}s\label{III-eq3.9}\nonumber\\
& & +I_{6}(t).\nonumber
\end{eqnarray}

Applying (\ref{II-eq2.5}) to $V(0,X^{\mu}_{0},X^{\nu}_{0})$, then
\begin{eqnarray}
\begin{split}
\frac{\|X^{\mu}_{t}-X^{\nu}_{t}\|^{2}}{\delta}
&\leq C\delta+\frac{C\|X^{\mu}_{0}-X^{\nu}_{0}\|^{2}}{\delta}+
 C\delta\int^{t}_{0}(1+L^{2}(s))\mathrm{d}s\label{III-eqq3.10}\\
&\ \ \ \  +C\int^{t}_{0}(1+L^{2}(s))(\frac{\|X^{\mu}_{s}-X^{\nu}_{s}\|^{2}}{\delta}+\frac{\mathcal{W}_{\overline{\mathcal{D}}_{s}}(\mu_{s},\nu_{s})^{2}}{\delta}\mathrm{d}s)\\
& \ \ \ \ +I_{6}(t).
\end{split}
\end{eqnarray}
By Doob's inequality,
\begin{eqnarray}\label{eq Zhai 5}
&&\ \ \ \ \ \mathbb{E}[\sup_{s\in[0,t]}\|I_{6}(s)\|^{2}]\\
&\leq&     16\int^{t}_{0}\mathbb{E}\left[\langle D_{x}u(s,X^{\mu}_{s},X^{\nu}_{s}),\sigma(s,X^{\mu}_{s},\mu_{s})\rangle ^{2}
        f^{2}_{\delta}(s,X^{\mu}_{s},X^{\nu}_{s})\right]\mathrm{d}s\nonumber\\
&& +16\int^{t}_{0}\mathbb{E}\left[\langle D_{y}u(s,X^{\mu}_{s},X^{\nu}_{s}),\sigma(s,X^{\nu}_{s},\nu_{s})\rangle ^{2}
        f^{2}_{\delta}(s,X^{\mu}_{s},X^{\nu}_{s})\right]\mathrm{d}s\nonumber\\
&& +16\int^{t}_{0}\mathbb{E}\left[\langle D_{x}f_{\delta}(s,X^{\mu}_{s},X^{\nu}_{s})
        +D_{y}f_{\delta}(s,X^{\mu}_{s},X^{\nu}_{s}),\sigma(s,X^{\mu}_{s},\mu_{s})\rangle ^{2}u^{2}(s,X^{\mu}_{s},X^{\nu}_{s})\right]\mathrm{d}s\nonumber\\
&& +16\int^{t}_{0}\mathbb{E}\left[\langle D_{y}f_{\delta}(s,X^{\mu}_{s},X^{\nu}_{s}),\sigma(s,X^{\mu}_{s},\nu_{s})
        -\sigma(s,X^{\nu}_{s},\mu_{s})\rangle ^{2}u^{2}(s,X^{\mu}_{s},X^{\nu}_{s})\right]\mathrm{d}s
        \nonumber\\
&\leq&  CT\delta^{2}\int^{t}_{0}(1+L(s))^{2}\mathrm{d}s+C\mathbb{E}\int^{t}_{0}(1+L(s))^{2}(\frac{\|X^{\mu}_{s}-X^{\nu}_{s}\|^{4}}{\delta^{2}}+
       \frac{\mathcal{W}_{\overline{\mathcal{D}}_{s}}(\mu_{s},\nu_{s})^{4}}{\delta^{2}})\mathrm{d}s.\nonumber
\end{eqnarray}

Letting $A_{t}:=\int^{t}_{0}(1+L^{2}(s))\mathrm{d}s$, by H\"older inequality, we have
$$(\int^{t}_{0}(1+L^{2}(s))f(s)ds)^{2}\leq A_{T}\int^{t}_{0}f^{2}(s)\mathrm{d}A_{s},\quad \forall f\in\mathcal{X}.$$
Squaring, multiplying $\delta^{2}$ and taking supremum and expectations on both sides of (\ref{III-eq3.9}), by the inequality above,
\begin{eqnarray}
\mathbb{E}[\sup_{s\in[0,t]}\|X^{\mu}_{s}-X^{\nu}_{s}\|^{4}]
&\leq& C\mathbb{E}[\|X^{\mu}_{0}-X^{\nu}_{0}\|^{4}]
+
C\delta^{4}+C\delta^{4}A^{2}_{T}\nonumber\\
& &+(C+CA_{T})\int^{t}_{0}\mathbb{E}[\sup_{u\in[0,s]}\|X^{\mu}_{u}-X^{\nu}_{u}\|^4\mathrm{d}]A_{s}\nonumber\\
& &+(C+CA_{T})\int^{t}_{0}\mathcal{W}_{\overline{\mathcal{D}}_{s}}(\mu_{s},\nu_{s})^{4}\mathrm{d}A_{s}\label{III-eq3.12}
\end{eqnarray}
Letting $\delta$ in \eqref{III-eq3.12} tending to zero, we have
\begin{eqnarray*}
\mathbb{E}[\sup_{s\in[0,t]}\|X^{\mu}_{s}-X^{\nu}_{s}\|^{4}]
&\leq& C\mathbb{E}[\|X^{\mu}_{0}-X^{\nu}_{0}\|^{4}]
+(C+CA_{T})\int^{t}_{0}\mathbb{E}[\sup_{u\in[0,s]}\|X^{\mu}_{u}-X^{\nu}_{u}\|^4]\mathrm{d}A_{s}\nonumber\\
& &+(C+CA_{T})\int^{t}_{0}\mathcal{W}_{\overline{\mathcal{D}}_{s}}(\mu_{s},\nu_{s})^{4}\mathrm{d}A_{s}\nonumber
\end{eqnarray*}
Then, by $X^{\mu}_{0}=X^{\nu}_{0}=X_0$, Gronwall's inequality and (\ref{III-eq3.1}), for any $t\in[0,T]$,
\begin{eqnarray}\label{eq Zhai 4}
\mathcal{W}_{t}(\Phi(\mu),\Phi(\nu))^{4}
&\leq&\mathbb{E}[\sup_{s\in[0,t]}\|X^{\mu}_{s}-X^{\nu}_{s}\|^{4}]\leq C\int^{t}_{0}\sup_{u\in[0,s]}\mathcal{W}_{\overline{\mathcal{D}}_{u}}(\mu_{u},\nu_{u})^{4}\mathrm{d}A_{s}\nonumber\\
&\leq&C\int^{t}_{0}\mathcal{W}_{s}(\mu,\nu)^{4}\mathrm{d}A_{s},
\end{eqnarray}
where $C$ is some constant dependent only on $T$.
Denote $\Phi^{k}$ as the $k$th composition of the mapping $\Phi$ with itself, for any $k\geq 1$,
\begin{eqnarray}
        \mathcal{W}_{T}(\Phi^{k}(\mu),\Phi^{k}(\nu))^{4}
&\leq&  C^{k}\int^{T}_{0}\frac{(A_{T}-A_{t})^{k-1}}{(k-1)!}\mathcal{W}_{t}(\mu,\nu)^{4}\mathrm{d}A_{t}\nonumber\\
&\leq&  \frac{(CA_{T})^{k}}{k!}\mathcal{W}_{T}(\mu,\nu)^{4}\label{III-eq3.14}.
\end{eqnarray}
For $k$ large enough, (\ref{III-eq3.14}) implies that $\Phi^{k}$ is strict contraction, there exists a unique fixed point for $\Phi^{k}$ on $\mathcal{P}^{2}(\mathcal{X})$. Thus there exists a unique strong solution $X$ to equation \eqref{I-eq1.2} on $[0,T]$.

  The proof is complete.
\end{proof}

\section{Propagation of Chaos}\label{III-sec3.2}
Let $n\in\mathbb{N}$ and denote $(X^{n,1}_{t},\cdots,X^{n,n}_{t})$ as a system of $n$ interacting particles driving by reflected SDEs with the form, for $i\in\{1,2,\cdots,n\}$,
\begin{eqnarray}
&&\mathrm{d}X^{n,i}_{t}=b(t,X^{n,i}_{t},\mu^{n}_{t})\mathrm{d}t+\sigma(t,X^{n,i}_{t},\mu^{n}_{t})\mathrm{d}W^{i}_{t}+\mathrm{d}K^{n,i}_{t},\ \mu^{n}_{t}=\frac{1}{n}\sum^{n}_{i=1}\delta_{X^{n,i}_{t}}\label{II-eq17},\\
& &X^{n,i}_{0}=x^{i},\quad |K^{n,i}|_{t}=\int_{0}^{t}\mathbf{1}_{\{X^{n,i}_{s}\in\partial \mathcal{D}_{s}\}}\mathrm{d}|K^{n,i}|_{s}, \quad K^{n,i}_{t}=\int_{0}^{t}\gamma(s,X^{n,i}_{s})\mathrm{d}|K^{n,i}|_{s},\nonumber
\end{eqnarray}
where the initial value $x^{1},\cdots, x^{n}$ are i.i.d. random variables with the same law as $X_0$ in (\ref{I-eq1.2}), $W^{1},W^{2},\cdots,W^{n}$ are mutually independent $m$-dimensional standard Wiener processes and $\delta_x$ is the Dirac measure concentrated at a point $x\in \mathbb{R}^d$.

\vskip 0.3cm

Now we demonstrate the propagation of chaos. Let $\mu^X$ be the law of the solution $X$ to equation \eqref{I-eq1.2} and $\mu^{n}:=\frac{1}{n}\sum^{n}_{i=1}\delta_{X^{n,i}}$ on $\mathcal{X}$.

\begin{thm}\label{II-thm 4.1}
Under Assumption (A1), the $n$-particle system given by (\ref{II-eq17}) converges in the following two sense. First,
\begin{equation}\label{II-eq2.18}
\lim_{n\to\infty}\mathbb{E}[\mathcal{W}_{\mathcal{X}}(\mu^{n},\mu^X)^{4}]=0.
\end{equation}
Second, for $n\in\mathbb{N}$, fixed $k\in\mathbb{N}$, the following weak convergence (or in distribution) holds, as $n\rightarrow\infty$,
\begin{equation}\label{II-eq2.19}
(X^{n,1},X^{n,2},\cdots,X^{n,k})\Rightarrow (Y^{1},Y^{2},\cdots,Y^{k}),
\end{equation}
where $Y^{1},\cdots,Y^{k}$ are the independent copies of the solution of (\ref{I-eq1.2}).
\end{thm}

\begin{proof}
Using the same Wiener processes and initial state $x^{i}$ as the $i$-th particle, we define $Y^{i}$ as the solution of the following reflected Mckean-Vlasov SDE
\begin{eqnarray*}
& & Y^{i}_{t}=x^{i}+\int^{t}_{0}b(s,Y^{i}_{s},\mu_{s})\mathrm{d}s+\int^{t}_{0}\sigma(s,Y^{i}_{s},\mu_{s})\mathrm{d}W^{i}_{s}+K^{i}_{t}, \quad \mu_{t}(dx)=\mathbb{P}[Y^{i}_{t}\in\mathrm{d}x], \nonumber
   \\
& & Y^{i}_{t}\in\overline{\mathcal{D}}_{t}, \quad |K^{i}|_{t}=\int^{t}_{0}\mathbf{1}_{\{Y^{i}_{s}\in\partial\mathcal{D}_{s}\}}\mathrm{d}|K^{i}|_{s}<\infty, \quad
    \ K^{i}_{t}=\int^{t}_{0}\gamma(s,Y^{i}_{s})\mathrm{d}|K^{i}|_{s}.
\end{eqnarray*}
It is easy to see that $\mu\equiv\mu^X$.
\vskip 0.3cm
Firstly, we obtain the estimation of the difference $\|X^{n,i}_{t}-Y^{i}_{t}\|^{4}$ given by
\begin{equation}\label{III-eq3.19}
        \mathbb{E}[\sup_{0\leq s\leq t}\|X^{n,i}_{s}-Y^{i}_{s}\|^{4}]
\leq C\mathbb{E}[\int^{t}_{0}\mathcal{W}_{s}(\mu^{n},\mu)^{4}\mathrm{d}A_{s}],
\end{equation}
due to the proof of Theorem \ref{II-thm2.2}; see (\ref{eq Zhai 4}).

Define $\nu^{n}:=\frac{1}{n}\sum^{n}_{i=1}\delta_{Y^{i}}$ as the empirical measure of $(Y^{i})_{ 0\leq i\leq n}$ and $\frac{1}{n}\sum^{n}_{i=1}\delta_{(X^{n,i},Y^{i})}$ is a coupling empirical measure of $\mu^{n}$ and $\nu^{n}$. Then
\begin{equation*}
\mathcal{W}_{t}(\mu^{n},\nu^{n})^{4}\leq\frac{1}{n}\sum^{n}_{i=1}\sup_{0\leq s\leq t}\|X^{n,i}_{s}-Y^{i}_{s}\|^{4}.
\end{equation*}
Combining this with (\ref{III-eq3.19}) to obtain
\begin{equation}
\mathbb{E}[\mathcal{W}_{t}(\mu^{n},\nu^{n})^{4}]\leq\mathbb{E}[\sup_{0\leq s\leq t}\|X^{n,i}_{s}-Y^{i}_{s}\|^{4}]\leq C\mathbb{E}[\int^{t}_{0}\mathcal{W}_{s}(\mu^{n},\mu)^{4}\mathrm{d}A_{s}].\label{III-eq3.20}
\end{equation}
According to the triangle inequality and (\ref{III-eq3.20}), we have
\begin{eqnarray}
        \mathbb{E}[\mathcal{W}_{t}(\mu^{n},\mu)^{4}]
&\leq&  8\mathbb{E}[\mathcal{W}_{t}(\mu^{n},\nu^{n})^{4}]
        +8\mathbb{E}[\mathcal{W}_{t}(\mu,\nu^{n})^{4}]
        \nonumber\\
&\leq& C\mathbb{E}[\int^{t}_{0}\mathcal{W}_{s}(\mu^{n},\mu)^{4}\mathrm{d}A_{s}]+C\mathbb{E}[\mathcal{W}_{t}(\mu,\nu^{n})^{4}].\nonumber
\end{eqnarray}
Then Gronwall's inequality implies that
\begin{eqnarray}
        \mathbb{E}[\mathcal{W}_{\mathcal{X}}(\mu^{n},\mu)^{4}] \leq C\mathbb{E}[\mathcal{W}_{\mathcal{X}}(\mu,\nu^{n})^{4}].\nonumber
\end{eqnarray}
Since $\nu^{n}$ is the empirical measures of i.i.d samples from the law $\mu$, we induce the limit in (\ref{II-eq2.18}) from the law of large numbers (see, e.g., \cite[Lemma 1.9]{Carmona}).
\vskip 0.3cm
Finally, by (\ref{III-eq3.19}) we can see that
\begin{eqnarray*}
        \mathbb{E}[\max_{1\leq i\leq k}\sup_{0\leq s\leq T}\|X^{n,i}_{s}-Y^{i}_{s}\|^{4}]
&\leq&  \sum^{k}_{i=1}\mathbb{E}[\sup_{0\leq s\leq T}\|X^{n,i}_{s}-Y^{i}_{s}\|^{4}]\\
&\leq&  C\mathbb{E}\left[\int^{T}_{0}\mathcal{W}_{t}(\mu^{n},\mu)^{4}\mathrm{d}A_{t}\right]\\
&\leq&  C\mathbb{E}[\mathcal{W}_{\mathcal{X}}(\mu^{n},\mu)^{4}].
\end{eqnarray*}
From the result in (\ref{II-eq2.18}), we know that $\mathbb{E}[\max_{1\leq i\leq k}\sup_{0\leq s\leq T}\|X^{n,i}_{s}-Y^{i}_{s}\|^{4}]$ converges to zero, then we get the claim limit in the (\ref{II-eq2.19}).
\end{proof}

\section{Large Deviation Principle}\label{section 5}
In this section, we assume that the initial data $X_0\in\overline{\mathcal{D}}_{0}$ is deterministic. For any $\varepsilon\in(0,1]$, let $X^{\varepsilon}=\{X^{\varepsilon}_{t},t\in[0,T]\}$ be
the unique strong solution to the following  reflected Mckean-Vlasov
SDE:
\begin{align}\label{IV-eq1}
\left\{
\begin{aligned}
& dX^{\varepsilon}_{t}=b(t,X^{\varepsilon}_{t},\mu^{\varepsilon}_{t})\mathrm{d}t+\sqrt{\varepsilon}\sigma(t,X^{\varepsilon}_{t},\mu^{\varepsilon}_{t})\mathrm{d}W_{t}
+K^{\varepsilon}_{t},\\
& X^{\varepsilon}_{0}=X_0,\quad \mu^{\varepsilon}_{t}(dx)=\mathbb{P}[X^{\varepsilon}_{t}\in\mathrm{d}x],\\
& |K^{\varepsilon}|_{t}=\int_{0}^{t}\mathbf{1}_{\{X^{\varepsilon}_{s}\in\partial \mathcal{D}_{s}\}}\mathrm{d}|K^{\varepsilon}|_{s},\quad  K^{\varepsilon}_{t}=\int_{0}^{t}\gamma(s,X^{\varepsilon}_{s})\mathrm{d}|K^{\varepsilon}|_{s}.
\end{aligned}
\right.
\end{align}
In this section, we mainly consider the LDP for $X^{\varepsilon}$ as $\varepsilon$ tending to $0$.

\vskip 0.3cm

We first recall the definition of LDP. Let $\{X^{\varepsilon}\}$ denote a family of random variables defined on $(\Omega,\mathcal{F},\mathbb{P})$ taking values in a Polish space $\mathcal{X}$.

\begin{dfn}(Rate function \cite{DZ,DupuisEllis})
A function $I: \mathcal{X}\to[0,\infty]$ is called a rate function on $\mathcal{X}$, if for any $C<\infty$, the level set $\{y\in \mathcal{X}: I(y)\leq C\}$ is a compact subset of $\mathcal{X}$.
\end{dfn}

\begin{dfn}(Large deviation \cite{DZ,DupuisEllis})
Let $I$ be a rate function on $\mathcal{X}$. The sequence $\{X^{\varepsilon}\}$ is said to satisfy LDP on $\mathcal{X}$ with the rate function $I$ if the following two conditions hold.

\begin{itemize}
  \item[(a)]{\bf Large Deviation upper bound.} For each closed subset $F$ of $\mathcal{X},$
  $$\limsup_{\varepsilon\to 0}\varepsilon\log\mathbb{P}\{X^{\varepsilon}\in F\}\leq -I(F).$$
  \item[(b)]{\bf Large Deviation lower bound.} For each open subset $G$ of $\mathcal{X},$
  $$\liminf_{\varepsilon\to 0}\varepsilon\log\mathbb{P}\{X^{\varepsilon}\in G\}\geq -I(G).$$
\end{itemize}
\end{dfn}

%
%
%
\vskip 0.3cm
From Theorem \ref{II-thm2.2}, we know that there exists a unique solution for the following reflected ODE
\begin{eqnarray}\label{IV-eq4.2}
& &\psi_{t}=X_0+\int^{t}_{0}b(s,\psi_{s},\delta_{\psi_{s}})\mathrm{d}s+K^{\psi}_{t},\quad  t\in [0,T],\\
& &|K^{\psi}|_{t}=\int_{0}^{t}\mathbf{1}_{\{\psi_{s}\in\partial \mathcal{D}_{s}\}}\mathrm{d}|K^{\psi}|_{s},\quad  K^{\psi}_{t}=\int_{0}^{t}\gamma(s,\psi_{s})\mathrm{d}|K^{\psi}|_{s}.\nonumber
\end{eqnarray}
We denote $\psi$ be the solution to (\ref{IV-eq4.2}).
For any $h\in L^{2}([0,T],\mathbb{R}^{m})$, consider the so called skeleton equation:
\begin{eqnarray}
& & Y^{h}_{t}=X_0+\int^{t}_{0}b(s,Y^{h}_{s},\delta_{\psi_{s}})\mathrm{d}s+\int^{t}_{0}\sigma(s,Y^{h}_{s},\delta_{\psi_{s}})h(s)\mathrm{d}s+K^{h}_{t},
\label{IV-eq4.3}\\
& & |K^{h}|_{t}=\int_{0}^{t}\mathbf{1}_{\{Y^{h}_{s}\in\partial \mathcal{D}_{s}\}}\mathrm{d}|K^{h}|_{s},\quad  K^{h}_{t}=\int_{0}^{t}\gamma(s,Y^{h}_{s})\mathrm{d}|K^{h}|_{s}.\nonumber
\end{eqnarray}
We stress that $\psi$ in \eqref{IV-eq4.3} is the strong solution to (\ref{IV-eq4.2}). Hence (\ref{IV-eq4.3}) is a classical, i.e. distribution independent, reflected  differential equation.

We have the following result:
\begin{prp}\label{V-prop5.1}
Under Assumption (A1), there exists a unique strong solution to equation (\ref{IV-eq4.3}).
\end{prp}
\begin{proof}
Let $\widetilde{b}(t,y):=b(t,y,\delta_{\psi_{t}})+\sigma(t,y,\delta_{\psi_{t}})h(t)$.
For any $t\in[0,T]$, $y,z\in \mathbb{R}^{d}$, by (\ref{II-eq2.16}),
\begin{eqnarray}
& &\|\widetilde{b}(t,y)-\widetilde{b}(t,z)\|\nonumber\\
&=&\|b(t,y,\delta_{\psi_{t}})-b(t,z,\delta_{\psi_{t}})+\sigma(t,y,\delta_{\psi_{t}})h(t)-\sigma(t,z,\delta_{\psi_{t}})h(t)\|\\
&\leq& L(t)(1+\|h(t)\|)\|y-z\|.\nonumber
\end{eqnarray}
Then, by Proposition \ref{prop 2} and the fact that $L(\cdot)(1+\|h(\cdot)\|)\in L^{1}([0,T],\mathbb{R}^+)$,  there exists a unique strong solution to (\ref{IV-eq4.3}).
\end{proof}

\vskip0.3cm

We now formulate the main result in this section as following Theorem.
\begin{thm}\label{II-thm5.2}
Let Assumption (A1) hold and $X^{\varepsilon}$ be the unique strong solution to (\ref{IV-eq1}). Then the family of $\{X^{\varepsilon}\}_{\varepsilon>0}$ satisfies a LDP on the space $\mathcal{X}$ with rate function
\begin{equation}
I(\phi):=\inf_{\{h\in L^{2}([0,T],\mathbb{R}^{m}):\phi=Y^{h}\}}\frac{1}{2}\int^{T}_{0}\|h(t)\|^{2}\mathrm{d}t,\ \phi\in\mathcal{X},
\end{equation}
with the convention $\inf\{\emptyset\}=+\infty$, here $Y^{h}\in\mathcal{X}$ solves equation (\ref{IV-eq4.3}).
\end{thm}
\begin{proof}
According to Proposition \ref{V-prop5.1}, there exists a measurable map
\begin{equation}\label{Gamma0}
\Gamma^{0}:C([0,T];\mathbb{R}^{m})\to \mathcal{X}\ \mbox{such\ that}\ Y^{h}=\Gamma^{0}(\int^{\cdot}_{0}h(s)\mathrm{d}s)\ \mbox{for}\ h\in L^{2}([0,T],\mathbb{R}^{m}).
\end{equation}
Let
$$\mathcal{H}^{M}:=\{h:[0,T]\to\mathbb{R}^{m}, \int^{T}_{0}\|h(s)\|^{2}\mathrm{d}s\leq M\},$$
and
\begin{equation}
\tilde{\mathcal{H}}^{M}:=\{h: h\ \mbox{is}\ \mathbb{R}^{m}\mbox{-valued}\ \mathcal{F}_{t}\mbox{-predictable process such that}\ h(\omega)\in \mathcal{H}^{M},\ \mathbb{P}\text{-}a.s\}.
\end{equation}
Throughout this section, $\mathcal{H}^{M}$ is endowed with the weak topology on $L^{2}([0,T],\mathbb{R}^{m})$. Then $\mathcal{H}^{M}$ is a Polish space.

By Theorem \ref{II-thm2.2} and \cite[Theorem 3.8]{LSZ} (or see (4.3) and (4.4) in \cite{LSZ}), for every $\varepsilon>0$, there exists a measurable mapping $\Gamma^{\varepsilon}(\cdot):C([0,T],\mathbb{R}^{m})\to\mathcal{X}$ such that  $$X^{\varepsilon}=\Gamma^{\varepsilon}(W(\cdot))$$
and  for any $M>0$ and $h^{\varepsilon}\in \tilde{\mathcal{H}}^{M}$,
$$Z^{\varepsilon}:=\Gamma^{\varepsilon}(W(\cdot)+\frac{1}{\sqrt{\varepsilon}}\int^{\cdot}_{0}h^{\varepsilon}(s)\mathrm{d}s)$$
is the solution of the following SDE
\begin{eqnarray}
    Z^{\varepsilon}_{t}
&=& X_0+\int^{t}_{0}b(s,Z^{\varepsilon}_{s},\mu^{\varepsilon}_{s})\mathrm{d}s
    +\sqrt{\varepsilon}\int^{t}_{0}\sigma(s,Z^{\varepsilon}_{s},\mu^{\varepsilon}_{s})\mathrm{d}W_{s}\nonumber\\
& & +\int^{t}_{0}\sigma(s,Z^{\varepsilon}_{s},\mu^{\varepsilon}_{s})h^{\varepsilon}(s)\mathrm{d}s+\int^{t}_{0}
    \gamma(s,Z^{\varepsilon}_{s})\mathrm{d}|K^{Z^{\varepsilon}}|_{s},\label{IV-eq4.17}\\
    |K^{Z^{\varepsilon}}|_{t}
&=& \int_{0}^{t}\mathbf{1}_{\{Z^{\varepsilon}_{s}\in\partial D_{s}\}}\mathrm{d}|K^{Z^{\varepsilon}}|_{s}, \quad K^{Z^{\varepsilon}}_{t}=\int_{0}^{t}
    \gamma(s,Z^{\varepsilon})\mathrm{d}|K^{Z^{\varepsilon}}|_{s}.\nonumber
\end{eqnarray}
Here $\mu^{\varepsilon}_{t}(dx)=\mathbb{P}[X^{\varepsilon}_{t}\in\mathrm{d}x]$, i.e.,  the law of $X^\epsilon$ on $\mathcal{X}$.

\vskip 0.3cm
According to Theorem 4.4 in \cite{LSZ}, to complete the proof of this theorem, it is sufficient to verify the following two claims:
\begin{itemize}
\item [(LDP1)]\label{LDP1} For every $M<+\infty$ and any family $\{h_{n};n\in\mathbb{N}\}\subset \mathcal{H}^{M}$
converging to some element $h$ in $\mathcal{H}^{M}$ as $n\to \infty$,
$$
\lim_{n\rightarrow\infty}\sup_{t\in[0,T]}\|\Gamma^{0}(\int^{\cdot}_{0}h_{n}(s)\mathrm{d}s)(t)-\Gamma^{0}(\int^{\cdot}_{0}h(s)\mathrm{d}s)(t)\|=0.
$$

\item [(LDP2)] For every $M<+\infty$ and any family $\{h^{\varepsilon};\varepsilon>0\}\subset \tilde{\mathcal{H}}^{M}$ and any $\theta>0$,
  $$\lim_{\varepsilon\to 0}\mathbb{P}(\sup_{t\in[0,T]}\|Z^{\varepsilon}_{t}-Y^{h^{\varepsilon}}_{t}\|>\theta)=0,$$
  where $Z^{\varepsilon}$ is the unique solution to (\ref{IV-eq4.17}) and $Y^{h^{\varepsilon}}=\Gamma^{0}(\int^{\cdot}_{0}h^{\varepsilon}(s)\mathrm{d}s)$.

\end{itemize}

The verification of (LDP1) will be given in Proposition \ref{IV-thm4.1}. (LDP2) will be established in Proposition \ref{IV-thm4.2}.
\end{proof}
\begin{remark}
We stress that $\mu^{\varepsilon}$ in \eqref{IV-eq4.17} is the distribution of strong solution to (\ref{IV-eq1}). This is somehow surprising.  The reason is that when perturbing the Brownian motion in the arguments of the mapping $\Gamma^{\varepsilon}(\cdot)$, $\mu^{\varepsilon}$ is already deterministic and hence it is not affected by the perturbation. For the details, we refer to Theorems 3.6, 3.8 and 4.4 in \cite{LSZ}.  An example is also introduced in \cite{LSZ};  see \cite[Example 1.1]{LSZ}.
\end{remark}

\subsection{Proof of (LDP1)}\label{IV-sub4.1}

In this subsection, we will prove the following result. 
\begin{prp}\label{IV-thm4.1}
Let Assumption (A1) hold. For any $M<+\infty$, family $\{h_{n}\}_{n\geq1}\subseteq\mathcal{H}^{M}$, and $h\in\mathcal{H}^{M}$, suppose that $\lim_{n\to\infty}h_{n}=h$ in the weak topology of $L^{2}([0,T];\mathbb{R}^{m})$. Then
\begin{equation}
\lim_{n\to\infty}\sup_{t\in[0,T]}\|\Gamma^{0}(\int^{\cdot}_{0}h_{n}(s)\mathrm{d}s)(t)-\Gamma^{0}(\int^{\cdot}_{0}h(s)\mathrm{d}s)(t)\|=0\label{IV-eq4.6}.
\end{equation}
\end{prp}

\begin{proof}
Let $(Y^{h},K^h)$ be the solution of (\ref{IV-eq4.3}), and $(Y^{h_{n}},K^{h_n})$ be the solution of (\ref{IV-eq4.3}) with $h$ replaced by $h_{n}$. By the definition of $\Gamma^{0}$, $Y^{h_{n}}=\Gamma^{0}(\int^{\cdot}_{0}h_{n}(s)\mathrm{d}s)$ and $Y^{h}=\Gamma^{0}(\int^{\cdot}_{0}h(s)\mathrm{d}s)$.
\vskip 0.3cm
The proof is divided into two steps.

 {\bf Step 1:} Set
        $\tilde{Z}^{h_{n}}_{t}=X_0+\int^{t}_{0}b(s,Y^{h_{n}}_{s},\delta_{\psi_{s}})\mathrm{d}s+\int^{t}_{0}\sigma(s,Y^{h_{n}}_{s},\delta_{\psi_{s}})h_{n}(s)\mathrm{d}s$,
  we claim that the set $\{\tilde{Z}^{h_{n}}\}_{n\geq 1}$ is pre-compact in $\mathcal{X}$.
       \vskip 0.3cm
        Due to the Cauchy-Schawrz inequality and (\ref{Remeq2.1}), there exists a constant $C_{M}$ depending on $M$ such that, for any $n\in\mathbb{N}$ and $0\leq s\leq t\leq T$,
        \begin{eqnarray}\label{IV-eq4.6'}
        \begin{split}
            &\ \ \ \ \|\tilde{Z}^{h_{n}}_{t}-\tilde{Z}^{h_{n}}_{s}\|\\
            &\leq  \|\int^{t}_{s}b(r,Y^{h_{n}}_{r},\delta_{\psi_{r}})\mathrm{d}r\|+
            \|\int^{t}_{s}\sigma(r,Y^{h_{n}}_{r},\delta_{\psi_{r}})h_{n}(r)\mathrm{d}r\|\\
            &\leq  \int^{t}_{s}\|b(r,Y^{h_{n}}_{r},\delta_{\psi_{r}})\|\mathrm{d}r +(\int^{t}_{s}\|\sigma(r,Y^{h_{n}}_{r},\delta_{\psi_{r}})\|^{2}\mathrm{d}r)
            ^{\frac{1}{2}}(\int^{t}_{s}\|h_{n}(r)\|^{2}\mathrm{d}r)^{\frac{1}{2}}\\
            &\leq  C_M(\int^{t}_{s}1+L^2(r)\mathrm{d}r)^{\frac{1}{2}}.
        \end{split}
        \end{eqnarray}
        To obtain the inequality above, $\{h_{n}\}_{n\geq1}\subseteq\mathcal{H}^{M}$ has been used.
        (\ref{IV-eq4.6'}), $L\in L^2([0,T],\mathbb{R}^+)$ and Remark \ref{Remark2.1} imply that the set $\{\tilde{Z}^{h_{n}}\}_{n\geq 1}$ is equicontinuous and $\sup_{n\in\mathbb{N}}\sup_{t\in[0,T]}\|\tilde{Z}^{h_{n}}_{t}\|<+\infty$.

Thus  $\{\tilde{Z}^{h_{n}}\}_{n\geq 1}$ is pre-compact in $\mathcal{X}$ due to the Arzela-Ascoli theorem.

\vskip 0.2cm
{\bf Step 2:} We verify that $\lim_{n\to\infty}\sup_{t\in[0,T]}\|Y^{h_{n}}_{t}-Y^{h}_{t}\|=0$, completing the  proof.
      \vskip 0.3cm
      By Lemma 4.4 in \cite{LO} and the result obtained in {\bf Step 1}, the set $\{Y^{h_{n}}\}_{n\geq1}$
      is relative compact in $\mathcal{X}$. Then there is a convergent subsequence of $\{Y^{h_{n}}\}_{n\geq 1}$,  which for notational convenience we again label
by $n$, and $\tilde{Y}\in\mathcal{X}$ such that
      \begin{eqnarray}\label{eq LDP0}
      \lim_{n\to\infty}\sup_{t\in[0,T]}\|Y^{h_{n}}_{t}-\tilde{Y}_{t}\|=0.
      \end{eqnarray}
      By Assumption (A1) and (\ref{eq LDP0}), for each $t\in [0,T]$,
        \begin{eqnarray}
            \lim_{n\to \infty}\int^{t}_{0}b(r,Y^{h_{n}}_{r},\delta_{\psi_{r}})\mathrm{d}r
            =\int^{t}_{0}b(r,\tilde{Y}_{r},\delta_{\psi_{r}})\mathrm{d}r.\label{IV-eq4.10}
        \end{eqnarray}
        For any $y\in\mathbb{R}^{d}$ and $t\in [0,T]$,
        \begin{eqnarray*}
        & &     |\langle \int^{t}_{0}\sigma(r,Y^{h_{n}}_{r},\delta_{\psi_{r}})h_{n}(r)\mathrm{d}r-\int^{t}_{0}\sigma(r,\tilde{Y}_{r},\delta_{\psi_{r}})h(r)\mathrm{d}r,y\rangle |\\
        &\leq&  (\int^{t}_{0}\|\sigma(r,Y^{h_{n}}_{r},\delta_{\psi_{r}})
        -\sigma(r,\tilde{Y}_{r},\delta_{\psi_{r}})\|^{2}
        \mathrm{d}r)^{\frac{1}{2}}(\int^{t}_{0}\|h_{n}(r)\|^{2}\mathrm{d}r)^{\frac{1}{2}}\|y\|\\
        & &     +|\langle\int^{t}_{0}\sigma(r,\tilde{Y}_{r},\delta_{\psi_{r}})(h_{n}(r)-h(r))\mathrm{d}r,y\rangle |\\
        &\leq&  \sqrt{M}|(\int_0^TL^2(t)dt)^{\frac{1}{2}}\|y\|\sup_{0\leq r\leq T}\|Y^{h_{n}}_{r}-\tilde{Y}_{r}\|\\
        & &     +|\langle \int^{t}_{0}\sigma(r,\tilde{Y}_{r},\delta_{\psi_{r}})(h_{n}(r)-h(r))\mathrm{d}r,y\rangle |.
        \end{eqnarray*}
        Combining this with the facts that $h_{n}$ converging to $h$ weekly in $L^2([0,T],\mathbb{R}^m)$ and $\langle\sigma(\cdot,\tilde{Y}_{\cdot},\delta_{\psi_{\cdot}}),y\rangle\in L^2([0,T],\mathbb{R}^m)$, and by Assumption (A1) and (\ref{eq LDP0}) again, we have
        \begin{equation}
        \lim_{n\to \infty}|\langle \int^{t}_{0}\sigma(r,Y^{h_{n}}_{r},\delta_{\psi_{r}})h_{n}(r)\mathrm{d}r
        -\int^{t}_{0}\sigma(r,\tilde{Y}_{r},\delta_{\psi_{r}})h(r)\mathrm{d}r,y\rangle |=0.\label{IV-eq4.12}
        \end{equation}
        Denote
        \begin{equation*}
        Z^{h}_{t}:=X_0+\int^{t}_{0}b(s,\tilde{Y}_{s},\delta_{\psi_{s}})\mathrm{d}s+\int^{t}_{0}\sigma(s,\tilde{Y}_{s},\delta_{\psi_{s}})h(s)\mathrm{d}s.
        \end{equation*}
        By (\ref{IV-eq4.10})-(\ref{IV-eq4.12}), for any $y\in\mathbb{R}^d$ and $t\in[0,T]$,  $\lim_{n\rightarrow\infty}\langle\tilde{Z}^{h_{n}}_{t},y\rangle=\langle Z^{h}_{t},y\rangle$. Combining this with
        the result obtained in {\bf Step 1}, we have
        \begin{eqnarray}\label{eq LDP 1}
        \lim_{n\rightarrow\infty}\sup_{t\in[0,T]}\|\tilde{Z}^{h_{n}}_{t}-Z^{h}_{t}\|=0.
        \end{eqnarray}

        By (\ref{eq LDP0}), (\ref{eq LDP 1}), and $K^{h_n}_{t}=Y^{h_{n}}_{t}-\tilde{Z}^{h_{n}}_{t},\ t\in[0,T]$, there is a $\tilde{K}\in C([0,T],\mathbb{R}^{d})$, such that $$\sup_{t\in[0,T]}|K^{h_{n}}_{t}-\tilde{K}_{t}|\to 0.$$
        Using a similar argument as in Lemma 4.5 in \cite{LO} and the uniqueness of the solution to (\ref{IV-eq4.3}), we have $(\tilde{Y},\tilde{K})=(Y^h,K^h)$. Then, by (\ref{eq LDP0}) again, $\lim_{n\to\infty}\sup_{t\in[0,T]}\|Y^{h_{n}}_{t}-Y^{h}_{t}\|=0$, completing the  proof.

%

\end{proof}

\subsection{Proof of (LDP2)}\label{IV-sub4.2}
For every $M<+\infty$ and any family $\{h^{\varepsilon};\varepsilon>0\}\subset \tilde{\mathcal{H}}^{M}$, recall that
$Z^{\varepsilon}:=\Gamma^{\varepsilon}(W(\cdot)+\frac{1}{\sqrt{\varepsilon}}\int^{\cdot}_{0}h^{\varepsilon}(s)\mathrm{d}s)$ and $Y^{h^{\varepsilon}}=\Gamma^{0}(\int^{\cdot}_{0}h^{\varepsilon}(s)\mathrm{d}s)$. In this subsection, we aim to prove the following result.
\begin{prp}\label{IV-thm4.2}
Let Assumption (A1) hold. Then for any $\theta>0$,
\begin{equation}
\lim_{\varepsilon\to 0}\mathbb{P}(\sup_{t\in[0,T]}\|Z^{\varepsilon}_{t}-Y^{h^{\varepsilon}}_{t}\|> \theta)=0.
\end{equation}
\end{prp}

Before proving Proposition \ref{IV-thm4.2}, we first prove the following a priori estimate.
\begin{lem}\label{IV-lem4.1}
Let $X^{\varepsilon}$ and $\psi$ satisfy (\ref{IV-eq1}) and (\ref{IV-eq4.2}), respectively. Then
\begin{eqnarray*}
 \lim_{\epsilon\rightarrow0}\sup_{t\in[0,T]}\mathcal{W}_{\overline{\mathcal{D}}_{t}}(\mu^{\varepsilon}_{t},\delta_{\psi_{t}})=0.
 \end{eqnarray*}
\end{lem}
\begin{proof}
Recall $V(t,x,y)$ introduced in (\ref{V}). Similar argument as proving (\ref{III-eqq3.10}) shows that, taking $\lambda$ large enough, for any $t\in[0,T]$,
\begin{eqnarray}
\begin{split}
 &\ \ \ \|X^{\varepsilon}_{t}
        -\psi_{t}\|^{2}\\
&\leq C\delta^2
         +
        C\int^{t}_{0}(1+L^2(s))(\|X^{\varepsilon}_{s}
        -\psi_{s}\|^{2}+\mathcal{W}_{\overline{\mathcal{D}}_{s}}(\mu^{\varepsilon}_{s},\delta_{\psi_{s}})^{2})\mathrm{d}s\\
        &\ \ \ +
        \delta J^\epsilon(t)
         +\delta\sqrt{\varepsilon}N^{\varepsilon}(t),
\end{split}
\end{eqnarray}
here
\begin{eqnarray*}
J^\epsilon(t)
=
\frac{1}{2}\varepsilon\int^{t}_{0}tr \left[\left(
                \begin{array}{cc}
                \sigma(s,X^{\varepsilon}_{s},\mu^{\varepsilon}_{s}) \\
                0 \\
                \end{array}
                \right)^{T}D^{2}V(s,X^{\varepsilon}_{s},\psi_{s})
                \left(
                \begin{array}{cc}
                \sigma(s,X^{\varepsilon}_{s},\mu^{\varepsilon}_{s}) \\
                0 \\
                \end{array}
                \right)\right]\mathrm{d}s,
\end{eqnarray*}
and
\begin{eqnarray*}
N^{\varepsilon}(t)  &:=&\int^{t}_{0}\langle D_{x}f_{\delta}(s,X^{\varepsilon}_{s},\psi_{s})
                        u(s,X^{\varepsilon}_{s},\psi_{s}), \sigma(s,X^{\varepsilon}_{s},\mu^{\varepsilon}_s)\rangle\mathrm{d}W_{s}\\
                    &&+\int^{t}_{0}\langle f_{\delta}(s,X^{\varepsilon}_{s},\psi_{s})D_{x}u(s,X^{\varepsilon}_{s},\psi_{s}),
                    \sigma(s,X^{\varepsilon}_{s},\mu^{\varepsilon}_s)\rangle\mathrm{d}W_{s}.
\end{eqnarray*}

The Gronwall lemma implies that
\begin{eqnarray}\label{eq LDP 4}
\begin{split}
 &\ \ \ \|X^{\varepsilon}_{t}
        -\psi_{t}\|^{2}\\
&\leq C\delta^2
         +
        C\int^{t}_{0}(1+L^2(s))\mathcal{W}_{\overline{\mathcal{D}}_{s}}(\mu^{\varepsilon}_{s},\delta_{\psi_{s}})^{2}\mathrm{d}s\\
        &\ \ \ +
       C \delta \sup_{s\in[0,T]}|J^\epsilon(s)|
         +C\delta\sqrt{\varepsilon}\sup_{s\in[0,T]}|N^{\varepsilon}(s)|.
\end{split}
\end{eqnarray}

Using an argument similar to proving (\ref{III-eq3.8}) and (\ref{eq Zhai 5}), we have
\begin{eqnarray*}
\sup_{s\in[0,T]}|J^\epsilon(s)|
\leq
C\epsilon\int_0^T(1+L^2(s))(\frac{1}{\delta}+\delta+\frac{\|X^{\varepsilon}_{s}-\psi_{s}\|^{2}}{\delta})ds,
\end{eqnarray*}
and
\begin{eqnarray*}
\mathbb{E}(\sup_{s\in[0,T]}|N^\epsilon(t)|^2)
\leq
C\int_0^T(1+L^2(s))(\delta^2+\frac{\mathbb{E}(\|X^{\varepsilon}_{s}-\psi_{s}\|^{4})}{\delta^2}+\frac{\mathbb{E}(\|X^{\varepsilon}_{s}-\psi_{s}\|^{2})}{\delta^2})ds.
\end{eqnarray*}
Combining the above two inequalities with the facts that $X^{\varepsilon}_{s}\in \overline{\mathcal{D}}_s$, $\psi_{s}\in \overline{\mathcal{D}}_s$, and $\mathcal{D}$ is a bounded domain in $\mathbb{R}^{1+d}$,
\begin{eqnarray}\label{eq LDP 2}
\sup_{s\in[0,T]}|J^\epsilon(s)|
\leq
C\epsilon(\frac{1}{\delta}+\delta),
\end{eqnarray}
and
\begin{eqnarray}\label{eq LDP 3}
\mathbb{E}(\sup_{s\in[0,T]}|N^\epsilon(s)|^2)
\leq
C(\delta^2+\frac{1}{\delta^2}).
\end{eqnarray}
By (\ref{eq LDP 4})--(\ref{eq LDP 3}), for any $t\in[0,T]$,
\begin{eqnarray}\label{eq LDP 5}
 \mathcal{W}_{\overline{\mathcal{D}}_{t}}(\mu^{\varepsilon}_{t},\delta_{\psi_{t}})^{2}
 &\leq& \mathbb{E}(\|X^{\varepsilon}_{t}
        -\psi_{t}\|^{2})\nonumber\\
&\leq& C\delta^2
         +
        C\int^{t}_{0}(1+L^2(s))\mathcal{W}_{\overline{\mathcal{D}}_{s}}(\mu^{\varepsilon}_{s},\delta_{\psi_{s}})^{2}\mathrm{d}s\nonumber\\
        && +
        C\delta \mathbb{E}(\sup_{s\in[0,T]}|J^\epsilon(s)|)
         +C\delta\sqrt{\varepsilon}\mathbb{E}(\sup_{s\in[0,T]}|N^{\varepsilon}(s)|)\nonumber\\
&\leq&
C\delta^2
         +
        C\int^{t}_{0}(1+L^2(s))\mathcal{W}_{\overline{\mathcal{D}}_{s}}(\mu^{\varepsilon}_{s},\delta_{\psi_{s}})^{2}\mathrm{d}s\nonumber\\
        && +
        C(1+\delta^2)\epsilon
         +C\sqrt{\varepsilon}(\delta^2+1).
\end{eqnarray}
Applying the Gronwall lemma, letting $\delta$ tend to $0$ and then $\varepsilon$ tend to $0$, we get
\begin{eqnarray}\label{eq LDP 6}
 \lim_{\epsilon\rightarrow0}\sup_{t\in[0,T]}\mathcal{W}_{\overline{\mathcal{D}}_{t}}(\mu^{\varepsilon}_{t},\delta_{\psi_{t}})^{2}=0.
 \end{eqnarray}

The proof is complete.

\end{proof}

Now we come back to the proof of Proposition \ref{IV-thm4.2}.
\begin{proof}
By Chebyshev's inequality, we only need to show that
\begin{equation}
\lim_{\varepsilon\to 0}\mathbb{E}\left[\sup_{t\in[0,T]}\|Z^{\varepsilon}_{t}-Y^{h^{\varepsilon}}_{t}\|^{2}\right]=0.
\end{equation}

By $\int^{T}_{0}\|h^{\varepsilon}(s)\|^2ds\leq M$ $\mathbb{P}$-a.s.,
\begin{eqnarray}\label{eq 2LDP 01}
&&\int^{T}_{0}(1+L(s))(1+\|h^{\varepsilon}(s)\|)ds\nonumber\\
&\leq& 2\int^{T}_{0}(1+L^2(s))ds+2\int^{T}_{0}(1+\|h^{\varepsilon}(s)\|^2)ds\leq C_M<\infty,\ \mathbb{P}\text{-a.s},
\end{eqnarray}
where $C_{M}$ is a constant depending on $M$.
By (\ref{eq 2LDP 01}), similar argument as proving (\ref{III-eqq3.10}) shows that, taking $\lambda$ large enough, for any $t\in[0,T]$,
\begin{eqnarray}\label{eq LDP00}
 &&\|Z^{\varepsilon}_{t}-Y^{h^{\varepsilon}}_{t}\|^{2}\nonumber\\
&\leq& C\delta^2
         +\nonumber
        C\int^{t}_{0}(1+L(s))(1+\|h^{\varepsilon}(s)\|)(\delta^2+\|Z^{\varepsilon}_{s}-Y^{h^{\varepsilon}}_{s}\|^{2}+\mathcal{W}_{\overline{\mathcal{D}}_{s}}(\mu^{\varepsilon}_{s},\delta_{\psi_{s}})^{2})\mathrm{d}s\\
        &&+
        \delta H^\epsilon(t)
         +\delta\sqrt{\varepsilon}A^{\varepsilon}(t)\nonumber\\
&\leq& C\delta^2
         +\nonumber
        C\int^{t}_{0}(1+L(s))(1+\|h^{\varepsilon}(s)\|)\|Z^{\varepsilon}_{s}-Y^{h^{\varepsilon}}_{s}\|^{2}\mathrm{d}s\\
        &&+C\sup_{s\in[0,T]}\mathcal{W}_{\overline{\mathcal{D}}_{s}}(\mu^{\varepsilon}_{s},\delta_{\psi_{s}})^{2}+
        \delta \sup_{s\in[0,T]}|H^\epsilon(s)|
         +\delta\sqrt{\varepsilon}\sup_{s\in[0,T]}|A^{\varepsilon}(s)|,
\end{eqnarray}
here
\begin{eqnarray*}
H^\epsilon(t)=\frac{1}{2}\varepsilon\int^{t}_{0}tr \left[\left(
                \begin{array}{cc}
                \sigma(s,Z^{\varepsilon}_{s},\mu^{\varepsilon}_{s}) \\
                0 \\
                \end{array}
                \right)^{T}D^{2}V(s,Z^{\varepsilon}_{s},Y^{h_{\varepsilon}}_{s})
                \left(
                \begin{array}{cc}
                \sigma(s,Z^{\varepsilon}_{s},\mu^{\varepsilon}_{s}) \\
                0 \\
                \end{array}
                \right)\right]\mathrm{d}s,
\end{eqnarray*}
and
\begin{eqnarray*}
      A^{\varepsilon}(t)
&=&    \int^{t}_{0}\langle D_{y}
        f_{\delta}(s,Z^{\varepsilon}_{s},Y^{h_{\varepsilon}}_{s})u(s,Z^{\varepsilon}_{s},Y^{h_{\varepsilon}}_{s}),\sigma(s,Z^{\varepsilon}_{s},
        \mu^{\varepsilon}_{s})\rangle\mathrm{d}W_{s}\\
& &     +\int^{t}_{0}\langle
        f_{\delta}(s,Z^{\varepsilon}_{s},Y^{h_{\varepsilon}}_{s})D_{y}u(s,Z^{\varepsilon}_{s},Y^{h_{\varepsilon}}_{s}),\sigma(s,Z^{\varepsilon}_{s},
        \mu^{\varepsilon}_{s})\rangle\mathrm{d}W_s.\nonumber
\end{eqnarray*}

Using similar arguments as proving (\ref{eq LDP 2}) and (\ref{eq LDP 3}),
\begin{eqnarray}\label{eq LDP 2-1}
\sup_{s\in[0,T]}|H^\epsilon(s)|
\leq
C\epsilon(\frac{1}{\delta}+\delta),
\end{eqnarray}
and
\begin{eqnarray}\label{eq LDP 3-1}
\mathbb{E}(\sup_{s\in[0,T]}|A^\epsilon(s)|^2)
\leq
C(\delta^2+\frac{1}{\delta^2}).
\end{eqnarray}
By (\ref{eq 2LDP 01})--(\ref{eq LDP 3-1}),
\begin{eqnarray*}
 \mathbb{E}(\sup_{t\in[0,T]}\|Z^{\varepsilon}_{t}-Y^{h^{\varepsilon}}_{t}\|^{2})
\leq
C\delta^2
         +
        C\sup_{s\in[0,T]}\mathcal{W}_{\overline{\mathcal{D}}_{s}}(\mu^{\varepsilon}_{s},\delta_{\psi_{s}})^{2}
        +
        C(1+\delta^2)\epsilon
         +C\sqrt{\varepsilon}(\delta^2+1).
\end{eqnarray*}

Applying Lemma \ref{IV-lem4.1}, letting $\delta$ tend to $0$ and then $\varepsilon$ tend to $0$, we get
$$
\lim_{\epsilon\rightarrow0}\mathbb{E}(\sup_{t\in[0,T]}\|Z^{\varepsilon}_{t}-Y^{h^{\varepsilon}}_{t}\|^{2})=0.
$$

The proof of Proposition \ref{IV-thm4.2} is complete.

\end{proof}

\vskip 0.4cm

\noindent{\bf  Acknowledgement.}\  This work is partly supported by the National Natural Science Foundation of China (No. 12131019, No. 11971456, No. 11671132, No. 11721101, No. 11871184) and the Fundamental Research Funds for the Central Universities (No. WK3470000024, No. WK0010000076).


\begin{thebibliography}{10}
\bibitem{AdamsReisRavaille}D.Adams, G.D.Reis, R.Ravaille, Large Deviation and Exit-times for reflected Mckean-Vlasov equations with self-stabilizing terms and superlinear drifts. Stochastic Processes and their Application, 2022, Volume 146, 264-310.

\bibitem{Bolley}F.Bolley, Separability and completeness for the Wasserstein distance.
Seminaire de probabilites XLI, Springer, 2008, pp.371-377.

\bibitem{BEH}P.Briand, R.Elie, Y.Hu, BSDEs with mean reflection. The Annals of Applied Probability, 2018, Volume 28, no.1, 482-510.

\bibitem{BCG}P.Briand, P.Chaudru de Raynal, A.Guillin, C.Labart, Particles systems and numerical schemes for mean reflected stochastic differential equations. The Annals of Applied Probability, 2020, Volume 30, no.4, 1884-1909.

\bibitem{BD2019}A.Budhiraja, P.Dupuis, Analysis and approximation of rare events: representations and weak convergence methods. Probability Theory and Stochastic Modeling, Springer, 2019, Volume 94.

\bibitem{BCS2004}K.Burdzy, Z.Chen, J.Sylvester, The heat equation and reflected Brownian motion in time dependent domains. The Annals of Probability, 2004, Volume 32, no.1B, 775-804.

\bibitem{[13]}K.Burdzy, D.Nualart, Brownian motion reflected on Brownian motion. Probability Theory and Related Fields, 2002, Volume 122, no.4, 471-493.

\bibitem{BurdzyW}K.Burdzy, W.N.Kang, K.Ramanan, The skorokhod problem in a time dependent interval.
Stochastic Processes and their Application, 2009, Volume 119, 428-452.

\bibitem{CDFM}M.Coghi, J.Deuschel, P.Friz, M.Maurelli, Pathwise McKean-Vlasov theory with additive noise. The Annals of Applied Probability, 2020, Volume 30, no.5, 2355-2392.

\bibitem{Costantini1990}C.Costantini, The Skorohod oblique reflection problem in
domains with corners and application to
stochastic differential equations. Probability Theory and Related Fields, 1992, Volume 91, 43-70.

\bibitem{Costantini2006}C. Costantini, E. Gobet, N. El Karoui,
Boundary sensitivities for diffusion processes in time dependent domains.
Applied Mathematics and Optimization, 2006, Volume 54, no.2, 159-187.

\bibitem{Carmona}R.Carmona, Lectures on BSDEs, Stochastic Control, and Stochastic Differential Games with Financial Applications. Society for Industrial and Applied Mathematics, 2016.

\bibitem{DZ}A.Dembo, O.Zeitouni. Large deviations techniques and applictions, Stochastic Modelling and Applied Probability. Springer-Verlag, Berlin 2010, Volume 38.

\bibitem{[RST]}G. Dos Reis, W. Salkeld, J. Tugaut, Freidlin-Wentzell LDPs in path space for
McKean-Vlasov equations and the functional iterated logarithm law. The Annals of Applied Probability, 2019, Volume 29, 1487-1540.

\bibitem{DupuisIshii}P.Dupuis, H.Ishii, SDEs with obliques reflection on nonsoomth domains, The Annals of probability, 1993, Volume 21, no.1, 554-580.

\bibitem{DupuisEllis}P.Dupuis, R.Ellis, A weak convergence approach to the theory of large deviations, Wiley, New York. 1997.

\bibitem{[27]}N. El Karoui, I. Karatzas, A new approach to the Skorohod problem, and its applications. Stochastics and Stochastics Reports, 1991, Volume 34, no.1-2, 57-82.

\bibitem{[28]}N. El Karoui, I. Karatzas, Correction: ``A new approach to the Skorohod problem, and its applications''.
Stochastics and Stochastics Reports, 1991, Volume 36, no.3-4, 265.

\bibitem{[26]}N. El Karoui, C. Kapoudjian, E. Pardoux, S. Peng, M.C. Quenez, Reflected solutions of backward SDE's, and related obstacle problems for PDE's. The Annays of Probability, 1997, Volume 25, no.2, 702-737.

\bibitem{GM1989}C. Graham, M. Metivier, System of interacting particles and nonlinear diffusion reflecting in a domain with sticky boundary. Probability Theory and Related Fields, 1989, Volume 82, no.2, 225-240.


\bibitem{HL}B. Hambly, S. Ledger,
A stochastic McKean-Vlasov equation for absorbing diffusions on the half-line.
The Annals of Applied Probability, 2017, Volume 27, no.5, 2698-2752.

\bibitem{[HIP]} S. Herrmann, P. Imkeller, D. Peithmann, Large deviations and a Kramers'
type law for self-stabilizing diffusions. The Annays of Applied Probabilty, 2008, 1379-1423.

\bibitem{[HLL2021]} W. Hong, S. Li, and W. Liu,  Large deviation principle for McKean-Vlasov quasilinear stochastic evolution equations. Applied Mathematicas and Optimization, 2021, Volume 84, no.1, suppl., S1119-S1147.

\bibitem{[36]} T. Konstantopoulos, V. Anantharam. An optimal flow control scheme that
regulates the burstiness of traffic subject to delay constraints. IEEE/ACM Transactions on
Networking, 1995, Volume 3, 423-432.
\bibitem{Lion}P.L.Lions, Stochastic Differential Equations with Reflecting Boundary Conditions. Communications on Pure and Applied Mathematics, 1984, Volume 37, 511-537.
\bibitem{LSZ}W.Liu, Y.L.Song, J.L.Zhai, T.S.Zhang, Large and Moderate Deviation vlasov principles for Mckean-vlasov SDEs with Jumps. To appear in Potential Analysis. Doi:10.1007/s11118-022-10005-0.
\bibitem{[SY]}Y. Suo and C. Yuan, Central Limit Theorem and Moderate Deviation Principle for McKean-Vlasov SDEs. Acta Applicandae Mathematicae, 2021, 175:16, 19 pp.
\bibitem{LO}N.L.P.Lundstr\"{o}m, T.\"{O}nskog, Stochastic and partial differential equations on non-smooth time-dependent domains. Stochastic Processes and their Applications, 2019, Volume 129, no.4, 1097-1131.
\bibitem{[43]} A. Mandelbaum, W.A. Massey, Strong approximations for time-dependent queues. Mathematics of Operations Research, 1995, Volume 20, no.1, 33-64.
\bibitem{MST}A.Matoussi, W.Sabbag, T.S.Zhang, Large deviation principle of obstacle problems for Quasilinear Stochastic PDEs. Applied Mathematics and Optimizations, 2021, Volume 83, no.2, 849-879.
\bibitem{Nystr"om}K. Nystr\"{o}m, T. \"{O}nskog, The Skorohod oblique reflection problem in time-dependent domains. The Annals of Probability, 2010, Volume 38, No.6, 2170-2223.

\bibitem{SaishoY}Y.Saisho, Stochastic Differential Equations for Multi-dimentional Domain with Reflecting Boundary. Probability Theory and Related Fields, 1987, Volume 74, 455-477.

\bibitem{Saisho}Y.Saisho, A model of the random motion of mutually reflecting molecules in $\mathbb{R}^{d}$, Kumamoto Journal of Mathematics, 1994, 95-123.

\bibitem{Skorohod1961}A.V. Skorohod, Stochastic equations for diffusion processes with a boundary. Theory of Probability \& Its Applications, 1961, Volume 6, 287-298.

\bibitem{Skorohod1962}A.V. Skorohod, Stochastic equations for diffusion processes with boundaries. II. Teory of Probability \& Its Applications, 1962, Volume 7, 5-25.

\bibitem{[59]}F. Soucaliuc, W. Werner, A note on reflecting Brownian motions. Electronic Communications in Probability, 2002, Volume 7, 117-122.

\bibitem{SznitamnA}A.A.Sznitamn, Nonlinear reflecting diffusion process, and the propagation of chaos and fluctuations associated.
Jounal of Functional analysis, 1984, Volume 56, no.3, 311-336.

\bibitem{Villani}C. Villani, Topics in optimal transportation, American Mathematical Society, 2003, no.58.

\bibitem{WZZ2022}R. Wang, J. Zhai, S. Zhang, Large deviation principle for stochastic Burgers type equation with reflection, 	
    Communications on Pure and Applied Analysis, 2022, Volume 21, 213-238.


\end{thebibliography}
\end{document}